\documentclass[reqno]{amsart}
\usepackage[english]{babel}
\usepackage{amsthm}
\usepackage{amsmath}
\usepackage{amssymb}
\usepackage[all]{xy}
\usepackage{dsfont}
\usepackage[bookmarksnumbered,colorlinks]{hyperref}
\usepackage{graphics}
\newcommand{\df}{{\rm d}}
\newcommand{\R}{\mathds R}

\def\ds{{\rm d}_s }
\def\dist{{\rm d} }

\def\d1#1#2{\frac{d#1}{d#2}}

\def\p1#1#2{\frac{\partial #1}{\partial #2}}
\def\to{t_o}

\def\part{a=\to \leq t_1 \leq \ldots \leq t_{n-1} \leq t_{n}=b}

\newcommand{\N}{\mathds N}

\newcommand\g{g} 
\newcommand\StStat{Stat}
\newcommand\GG{M}

\title[Interplay between Lorentzian Causality and Finsler metrics]{On the interplay between Lorentzian\\ Causality and  Finsler metrics \\ of Randers type}
\author[E. Caponio]{Erasmo Caponio}
\address{Dipartimento di Matematica, Politecnico di Bari, Via Orabona 4,
70125, Bari, Italy}
\email{caponio@poliba.it}
\thanks{E.C. supported by M.I.U.R. Research project PRIN07 ``Metodi Variazionali e Topologici nello Studio di Fenomeni Nonlineari''
}
\author[M. A. Javaloyes]{Miguel Angel Javaloyes}
\address{Departamento de Geometr\'{\i}a y Topolog\'{\i}a.
 Facultad de Ciencias, Universidad de Granada.
 Campus Fuentenueva s/n, 18071 Granada, Spain}
\email{ma.javaloyes@gmail.com}
\curraddr{Departamento de Matem\'aticas, Facultad de Matem\'aticas, Universidad de Murcia,
Campus Universitario de Espinardo, 30100 Murcia, Spain}
\thanks{MAJ is partially supported by Regional J. Andaluc\'{\i}a Grant P09-FQM-4496, MICINN
project MTM2009-10418, and Fundaci\'on S\'eneca project 04540/GERM/06. This research is a result of the
activity developed within the framework of the Programme in Support of Excellence Groups of the Regi\'on de Murcia, Spain, by Fundaci\'on S\'eneca, Regional Agency for Science and Technology (Regional
Plan for Science and Technology 2007–2010).}

\author[M. S\'anchez]{Miguel S\'anchez}
\address{Departamento de Geometr\'{\i}a y Topolog\'{\i}a.
 Facultad de Ciencias, Universidad de Granada.
 Campus Fuentenueva s/n, 18071 Granada, Spain}
\email{sanchezm@ugr.es}

\thanks{MS is partially supported by
Spanish MEC-FEDER Grant MTM2007-60731 and Regional J. Andaluc\'{\i}a Grant P09-FQM-4496}

\thanks{All the authors are partially supported by the Spanish-Italian Acci\'on Integrada HI2008.0106/Azione Integrata
Italia-Spagna IT09L719F1}



\date{}

\begin{document}
\newtheorem{thm}{Theorem}[section]
\newtheorem{prop}[thm]{Proposition}
\newtheorem{lemma}[thm]{Lemma}
\newtheorem{cor}[thm]{Corollary}
\theoremstyle{definition}
\newtheorem{defini}[thm]{Definition}
\newtheorem{notation}[thm]{Notation}
\newtheorem{exe}[thm]{Example}
\newtheorem{conj}[thm]{Conjecture}
\newtheorem{prob}[thm]{Problem}
\newtheorem{rem}[thm]{Remark}

\begin{abstract}
We obtain some results in both Lorentz and Finsler geometries, by
using a correspondence between the conformal structure
(Causality) of standard stationary spacetimes on $M=\R\times S$
and Randers metrics on $S$. In particular:

(1) For stationary spacetimes: we give a simple characterization
of when $\R\times S$ is causally continuous or globally hyperbolic
(including in the latter case, when $S$ is a Cauchy hypersurface),
in terms of an associated Randers metric. Consequences for the
computability of Cauchy  developments  are also derived.

(2) For Finsler geometry:  Causality suggests that the role of
completeness in many results of Riemannian Geometry (geodesic
connectedness by minimizing geodesics, Bonnet-Myers, Synge
theorems) is played by the compactness of symmetrized closed balls in
Finslerian Geometry.
Moreover, under this condition we show that for any Randers metric
$R$
there exists another
Randers metric $\tilde R$ with the same pregeodesics and
geodesically complete.

Even more, results on the differentiability of Cauchy horizons in
spacetimes yield consequences for the differentiability of the
Randers distance to a subset, and vice versa.

\end{abstract}

\maketitle
\section{Introduction}

Randers  metrics constitute the most typical class of
non-reversible Finsler metrics, and the differences between their
properties and those of Riemannian metrics become apparent. For
example, they include compact Katok manifolds, which admit only
finitely many closed geodesics (see \cite{ka73,Zi83}); in
particular, there are Katok metrics on the sphere $S^2$ with only
two distinct closed geodesics, whereas any Riemannian metric on
$S^2$ admits infinitely many  (see \cite{Ban93,Franks92}). There
are several ways to express a given Randers metric; for instance,
 when considered as a Zermelo metric, the characterization of those
with constant flag curvature becomes more natural (see
\cite{BCS04}). Moreover, there is an interesting relation between
standard stationary spacetimes $(\R\times S,\g)$ (see Eq.
\eqref{standard}) and Randers metrics (see Eq. \eqref{fermat}).
This was pointed out in \cite{cmm}, where the associated Randers
metric is called   {\em Fermat metric}, and was used in this
reference and others \cite{BaCaCa08, BiJa08, cmm2} to prove some
properties of geodesics in stationary spacetimes. Such Fermat
metrics are  referred to  as Optical Zermelo-Randers-Finsler
metrics  in \cite{GHWW}, where some interpretations (concerning,
for example, the case of constant flag curvature) are provided.

The purpose of the present paper is to develop in full the
correspondence between the global conformal properties (Causality)
of the stationary metric and the global geometric properties of
the associated Fermat/optical metric. This will be useful for both
geometries, Lorentzian and Finslerian, and it is summarized now.



\subsection{Previous notions on stationary spacetimes} For the
convenience of the reader,  we recall first some basic elements on
stationary spacetimes which are necessary in order to understand
the announced correspondence with Randers metrics. Such elements
are spread in physics and mathematics literature (see for example
\cite{CFS, FS,Mas}, \cite[pp. 37-38]{oneill95}, \cite[Ch.
18]{SKMHH}) and some of them will be developed in more detail in
Section \ref{stationaryfermat} (see also Section \ref{causality}
for general background on Lorentzian Geometry).

A (normalized, standard) stationary spacetime is a  smooth, connected, product
manifold $ M=\R\times S$ endowed with a Lorentzian metric $g$
(with signature $(-,+,\dots,+)$) which can be written as:
\begin{equation}\label{e1} g=-dt^2+ \pi^*\omega \otimes dt + dt \otimes \pi^*\omega  +\pi^*g_0 ,
\end{equation}
where $\pi: \R\times S \rightarrow S$, $ t: \R\times S
\rightarrow \R$ are the natural projections, $^*$ denotes
pullback and $\omega$ and $g_0$ are resp. a $1$-form and a Riemannian metric on $S$.
Here, the vector field $K$ induced from the natural
lifting to $M$ of the canonical  vector field on  $\R$ is a
(normalized, standard) timelike Killing vector field.

Let {\rm \StStat}$(\R\times S)$ be the set of all the standard
stationary metrics on $\R\times S$ as in (\ref{e1}). Notice that
any such metric is determined by the pair  $(g_0,\omega)$
composed by a Riemannian metric $g_0$ and a $1$-form $\omega$ on
$S$; such  a pair will be called a {\em stationary data pair}.
Conversely, any stationary data pair $(g_0,\omega)$ determines a
Lorentzian metric through the expression (\ref{e1}) and, so, it
defines an element of {\rm \StStat}$(\R\times S)$. Therefore, {\rm
\StStat}$(\R\times S)$ can be also regarded as the set of all the
stationary data pairs  on $S$. We emphasize that two distinct stationary
data  pairs
may yield isometric spacetimes. This is straightforward because an
expression such as (\ref{e1}) can be derived by making any
spacelike section $S'$ to play the role of the initial section
$S$. Nevertheless, we detail this assertion for clarity in future
reference. Start with the metric (\ref{e1}) obtained for some
stationary data pair $(g_0,\omega)$. On any embedded hypersurface
$i:S'\hookrightarrow M$, one can induce a (possibly
signature-changing) metric $g_0'=i^*g$, and $S'$ is called {\em
spacelike} if  $g_0'$ is a Riemannian (positive-definite) metric.
Moreover, a $1$-form $\omega'$ can be induced as $\omega'(v)=g(K,v)$
for all $v$ tangent to $S'$. Assume that $S'$ is spacelike and
also a section, that is, $S'$ can be written as the
graph
 $S^f=\{(f(x),x):
x\in S\}$ of some smooth function $f$ on $S$. The restriction
$\pi|_{S^f}: S^f\rightarrow S$ is a diffeomorphism, and we write
$\tilde f=(\pi|_{S^f})^{-1}$, 
$g^f_0= \tilde f^*g'_0$ and 
$\omega^f= \tilde f^*\omega'$. 
One can check that the $1$-forms  vary always in the same cohomology
class as $\omega - \omega^f=df$, and the metrics satisfy $g_0 +
\omega\otimes\omega=g_0^f + \omega^f\otimes\omega^f$ (Prop.
\ref{multfermat}).
Now, the
metric $g^f$ in {\rm \StStat}$(\R\times S)$ determined by
 the stationary data pair $(g_0^f,\omega^f)$ 
is isometric to the original one. In fact, the {\em change of
initial section} $f_M: M\rightarrow M$
\begin{equation}\label{F}
f_M(t,x)=(t+f(x),x) \quad \quad \forall (t,x)\in \R\times S
\end{equation}
satisfies
$f_M^*g=g^f$.

\subsection{Stationary-Randers Correspondence (SRC)}\label{ssrc} Notice that any stationary data pair
$(g_0,\omega)$   also determines a Finsler metric of Randers type
on S, namely:
$$  R(v)= \sqrt{h(v,v)} + \omega(v), \quad \quad \hbox{for } \;
 h=g_0+\omega\otimes \omega,\ v\in TS.$$ Any Randers metric on $S$ can be obtained in
such a way. Therefore, the set of all such pairs $(g_0,\omega)$
can be also identified with the set  {\rm Rand}$(S)$ of all the
Randers metrics on $S$, and one has the natural bijective map
\begin{equation}\label{eee}
{\rm \StStat} (\R\times S)  \rightarrow  {\rm Rand}(S), \quad
\quad g \mapsto F_{\g}, \end{equation} where  $F_g$ is determined
by the same stationary data pair $(g_0,\omega_0)$ which determines
$g$. We call $F_g$ the {\em Fermat metric} associated to $g$.

Now, some first relations  between the geometric properties of the
spacetime $(M,g)$ and the Randers manifold $(S,F_g)$ appear, for
example:
\begin{itemize}\item[(A1)] The future-pointing (resp. past-pointing) lightlike
pregeodesics of $(M,g)$ project onto the pregeodesics (resp.
reverse pregeodesics) of $(S,F_g)$ (Prop. \ref{ppp}). \item[(A2)]
A slice $S_t=\{t\}\times S$ is a Cauchy hypersurface for $(M,\g)$
(and, then, $(M,  \g)$ is globally hyperbolic) if and only if the
corresponding Fermat metric $F_{\g}$ is forward and backward
complete (Th. \ref{th2}).
\end{itemize}
We can go further in this correspondence by noticing the
following claim.
 Assume that two stationary data pairs $(g_0,\omega),
(g'_0,\omega')$ yield stationary metrics $g$, $g'$ which are
isometric by means of a change of the initial section as in
(\ref{F}). {\em Then, the associated Fermat metrics $F_g$,
$F_{g'}$
 must share the
geometric properties which correspond to the intrinsic properties
of $g$.}
 More precisely,
define the following relations of equivalence in {\rm Rand}$(S)$
and {\rm \StStat}$(\R\times S)$, resp.:
$$\begin{array}{clll}
R\sim R'& \Longleftrightarrow &  R-R'= df &\hbox{for some smooth function $f$ on $S$},\\
\g\sim \g ' & \Longleftrightarrow & g'=f_M^*g &\hbox{for some
change of the initial section $f_M$}\\ 
&&&\hbox{as in (\ref{F}),}
\end{array}$$
and let  {\rm Rand}$(S)/\sim$, {\rm \StStat}$(\R\times S)/\sim$ be
the corresponding quotient sets, resp. The bijection (\ref{eee})
induces a well-defined bijective map between the quotients
$$({\rm \StStat}(\R\times S)/\sim ) \rightarrow {\rm
(Rand}(S)/\sim )
$$
 (see Prop. \ref{ppp} and \ref{multfermat}). Then,  the claim above yields, for example
(see Th. \ref{causalstructure}):

\begin{itemize}
\item[(B1)] One (and then all)  representative of the {\em class}  $[\g]$ $(\in {\rm
\StStat}(\R\times S)/\sim )$ is globally hyperbolic if and only if
 one (and then  all) representative of $[F_{\g}] (\in
{\rm Rand}(S)/\sim )$ satisfies that  its symmetrized closed balls  (i.e. the closed balls defined by the symmetrized distance associated to $F_{\g}$)  are compact. In particular, this happens if the Randers distance is forward or
backward complete for some representative $F_{\g}\in [F_{\g}]$.

\item[(B2)] One (and then all) representative of $[\g]$ is
causally simple if and only if one (and then all) representative
of $[F_{\g}]$ is convex,  i.e. any two points in $S$ can be joined by a geodesic whose length is equal to the distance between the two points. As a consequence, $\g$ satisfies that
$J^+(p)$ is closed for all $p$ if and only if $J^-(p)$ is closed
for all $p$.
\end{itemize}

\subsection{Applications to standard stationary spacetimes.}
SRC yields explicit applications for the study of the so-called
``causal ladder of spacetimes'', and this can be extended to other
causal properties. Recall that standard stationary spacetimes are
always causally continuous, and   the next conditions in  the
causal ladder are
 causal simplicity and global hyperbolicity (see Section \ref{causality}).
 In general, these conditions may be difficult to check, but items (B1) and (B2) yield a complete
characterization for standard stationary spacetimes\footnote{As
these properties are conformally invariant, the normalization
$g(K,K)=-1$ for $\g$ in SRC is just a non-restrictive choice.
Along the remainder of the paper, the results will be written
without this normalization in order to make the expressions valid
directly for all stationary spacetimes (and also trivially for all
conformastationary spacetimes).}. We emphasize that,  in (B1), global
hyperbolicity is characterized even when $S_0= \{0\}\times S$ is not a Cauchy
hypersurface of $\R\times S$ (typically, global hyperbolicity is
proved by checking that some candidate hypersurface is Cauchy but,
in principle, the only candidate in our case would be $S_0$ or, equivalently, any slice
$S_t= \{t\}\times S, t\in \R$).
However, one can check directly when $S_0$ is a Cauchy hypersurface
through the characterization in item (A2). This sharpens the
natural rough estimates for Cauchy hypersurfaces in splitting type
spacetimes, see \cite{Sa-bari}.

An example of other causal properties  which can be studied via
Finslerian ones will be developed in Subsection \ref{s4.3}. Here a
simple application for the computation of  Cauchy developments
(Prop. \ref{p4.5}), including the problem  of the (non-)smoothness
of the Cauchy horizon (Th. \ref{nondiffmeasure}), is obtained.

\subsection{Applications to Finsler geometry.} They appear in
many directions. First, the plain applications to
Causality work also the other way round; so, results on, say,
differentiability of Cauchy horizons, are translated in results on
the differentiability of the distance function to a set for a
Randers metric (Th. \ref{nondiff} and Cor. \ref{nullmeasure}).

 Causality also suggests the appropriate hypotheses to study
the geometry of Randers metrics and, eventually, for any Finsler
metric. Indeed,  from SRC an analogy between Riemannian and
Randers metrics for the problem of convexity  becomes clear: {\em
the role of the compactness of the symmetrized closed balls for a
Randers metric is  similar to the role of metric completeness for
a Riemannian metric}. In fact, as globally hyperbolic spacetimes
are causally continuous, the items (B1) and (B2) of SRC imply
directly that any Randers metric with compact symmetrized closed
balls is convex (see Lemma \ref{lgeodesicconnection} for details).
Moreover, this property can  be also generalized to any Finsler
metric, as it is carried out by using variational arguments in Th.
\ref{tgeodesicconnection}. Analogously, one can extend the
Finslerian theorem of Bonnet-Myers (see \cite[Theorems
7.7.1]{BCS01}) as well as other theorems where this is used (for
example, Synge or the  sphere theorem). Namely, forward or
backward completeness can be replaced by the weaker condition of
compactness for the symmetrized closed balls.

A different type of applications of SRC appears for those Randers
metrics $ R(v)=\sqrt{h(v,v)}+\omega(v)$ with equal $h$ but different $\omega$
in the same cohomology class. As claimed in Subsection
\ref{ssrc}, all of them are Fermat metrics for canonically
isometric
 standard stationary spacetimes and, therefore,
 they share the elements which correspond
to the (intrinsic) geometry of these spacetimes. The first of
these elements is the set of pregeodesics (which corresponds to
the lightlike pregeodesics of the spacetime, item (A1)).  Another
shared property is the possible compactness of symmetrized closed
balls (as this corresponds to global hyperbolicity, item (B1)) as
well as convexity (as this corresponds to causal  simplicity, item
(B2)).

More subtly,  recall that the fact that each slice $ S_t=\{t\}\times S$ is a
Cauchy hypersurface, is not an intrinsic property for the spacetime
(say, it depends on the
 choice of the initial spacelike hypersurface $S_0$). In spite of this, such a  property was also
characterized in terms of the completeness of the corresponding
Fermat metric (item (A2)).  Therefore, one obtains a surprising
consequence for Randers metrics: {\em the class of a Randers
metric $R$  (according to Subsection \ref{ssrc}) contains a forward and backward complete
representative (necessarily with the same pregeodesics  as $R$) if
and only if its  symmetrized closed balls are compact},
see\footnote{In comparison, notice that a well-known result by
Nomizu and Ozeki \cite{NO} states that any Riemannian manifold is
globally conformal to a complete one --with different
pregeodesics, in general.  If a Riemannian metric is incomplete
then, regarded as a Finsler metric, it contains non-compact
(symmetrized) closed
 balls and, so,  no Randers metric in its class will
be complete by SRC (B1).} Th. \ref{randerscompleteness}. We recall
that, even though the symmetrized distance $\ds$ appears sometimes
in the literature (see for example \cite{Rad04,Rad04b}), there are
no similar results to previous ones,  as far as we know. A
difficulty of $\ds$ is that it is not constructed as the distance
associated to a length structure (see Appendix).

 The interaction is also symbiotic in other
intermediate aspects. For example, under  the Finslerian
viewpoint, the completeness of the symmetrized distance $\ds$ of
the Randers metric $ R(v)=\sqrt{h(v,v)}+\omega(v)$  is easily a sufficient
(but not necessary) condition for the completeness of the
associated Riemannian metric $h$. This yields that, if  a standard
stationary spacetime  $(M,g)$ is globally hyperbolic, then
 the Riemannian metric  $g_0+ \omega\otimes \omega$
  on $S$,
(the ``quotient metric by the flow of $K$'', under the normalization
$g(K,K)=-1$), which is the Riemannian metric $h$ of the Fermat metric $F_g$, must be complete (see Eq. \eqref{standard} and Cor.
\ref{globhyp}), and counterexamples to the converse
   follow easily, Example \ref{nosymcomplete}.
\subsection{Plan of work} Due  to the quite interdisciplinary nature of this work,
a special effort  has been carried out in order to keep it
self-contained. The paper is organized as follows. In Section
\ref{finslermetrics} we introduce the basic notions and results in
Finsler manifolds necessary for the remainder of the paper.
Moreover, we give an extension of the Hopf-Rinow theorem for
Finsler manifolds which considers the symmetrized balls  (see
Prop. \ref{dscomplete}). In Section \ref{causality}, we recall the
basic notions on Causality in spacetimes, including the causal
ladder. In Section \ref{stationaryfermat},  the Fermat metric
associated to any (non-necessarily normalized) standard stationary
spacetime is introduced (see \eqref{fermat} and \eqref{standard}).
Then, the applications to Causality of standard stationary
spacetimes are obtained (Th. \ref{causalstructure}, \ref{th2}).
Moreover, in Prop. \ref{p4.5} we characterize the future/past
Cauchy development of a subset $A$ contained in a Cauchy
hypersurface, by using the distance function
 from $A$ in the Fermat metric. As a consequence, we give a
result about the measure of the subset of  non-differentiable
points in the Cauchy horizons $H^\pm(A)$ when $A$ is a domain with
enough regular edge (Th. \ref{nondiffmeasure}). In Section
\ref{randersfermat}, we exploit the expression of Randers metrics
as Fermat ones in order to obtain the applications to Finsler
Geometry. First, we prove the existence and multiplicity of
connecting geodesics for general Finsler manifolds (Th.
\ref{tgeodesicconnection}) taking into account the hypothesis on
Randers metrics suggested by Causality (Lemma
\ref{lgeodesicconnection}). Second, we retrieve a necessary
condition for a standard stationary spacetime to be globally
hyperbolic (Cor. \ref{globhyp}), and obtain a result on geodesic
completeness for Randers metrics with compact symmetrized closed
balls (Th. \ref{randerscompleteness}).
We finish the section by applying
some results on differentiability of Cauchy horizons in
\cite{CFGH} to the differentiability of the distance function to a
subset in a Randers metric (Th. \ref{nondiff} and Cor.
\ref{nullmeasure}). In  the Appendix, we discuss the relation
between the different metrics which appear in Randers manifolds,
introducing the length metric  associated to the symmetrized
distance.

As possible further developments, notice that most of the results
in this paper apply only to Randers metrics. It is natural to
wonder which of them can be extended to arbitrary (non-reversible)
Finsler manifolds.

\section{Finsler metrics}\label{finslermetrics}
Let $M$ be a $ C^\infty$, paracompact, connected
manifold of dimension $n$ and $F:TM\rightarrow [0,+\infty]$  be a
continuous function. We say that $(M,F)$ is a Finsler manifold if
\begin{enumerate}
\item $F$ is $C^\infty$ in $TM\setminus 0$, i. e. away from the
zero section, \item $F$ is fiberwise positively homogeneous of
degree one, i. e. $F(x, \lambda y)=\lambda F(x, y)$ for every
$(x,y)\in TM$ and $\lambda>0$, \item $F^2$ is fiberwise strongly
convex, i. e. the matrix
\begin{equation}\label{fundamentaltensor}
g_{ij}(x,y)=\left[\frac{1}{2}\frac{\partial^2 (F^2)}{\partial
y^i\partial y^j}(x,y)\right]
\end{equation}
is positive definite for every $(x,y)\in TM\setminus 0$.
\end{enumerate}
Here we are using the  notation of \cite{BCS01}, that is, $(x,y)$
denotes the natural coordinates in $TM$ associated to a chart in
$M$. In what follows, $v$ will denote directly a vector in $TM$
and the reference to the base point (which eventually will be
denoted with a different letter, $p, q \in M$) will be omitted
when there is no possibility of confusion.

Given a Finsler manifold $(M,F)$, it can be proven that
$F$ must be in fact positive away from the zero section, and it
satisfies the {\it triangle inequality}, that is,
\begin{equation}\label{triineq}
F(v_1+v_2)\leq F(v_1)+F(v_2),
\end{equation}
where equality holds iff $v_2=\alpha v_1$ or $v_1=\alpha v_2$ for some $\alpha\geq 0$, and the {\it fundamental inequality},
\[
\sum_{ij} g_{ij}(v)w^iv^j\leq F(w)F(v),
\]
where $v\not=0$, and equality holds {iff} $w=\alpha v$ for some $\alpha\geq 0$.

A remarkable property of Finsler metrics is that they may be
non reversible, that is, in general $F(-v)\not=F(v)$. So, one
defines the {\it reverse Finsler metric} $\tilde{F}$  by
$\tilde{F}(v)=F(-v)$ for every $v\in TM$, which is again a
Finsler metric. The most typical non-reversible examples are
Randers metrics. Let  $(M,h)$ be a Riemannian manifold and
$\omega$ be a $1$-form on $M$  such that $|\omega(v)|<\sqrt{h(v,v)}$ for any $v\in TM$. The Randers metric $R$ on $(M,h)$ associated to
$\omega$ is then
\begin{equation}\label{randers}
 R(v)=\sqrt{h(v,v)}+\omega(v),
\end{equation}
(see \cite[Chapter 11]{BCS01}). The reverse metric $\tilde R$ is
obtained just replacing $\omega$ by $-\omega$.

\subsection{Distance function and length.} For any Finsler metric, one can define
naturally the (Finslerian) distance as
\begin{equation}\label{distance}
\dist (p,q)=\inf_{\gamma\in C(p,q)} \ell_F(\gamma),
\end{equation}
where $C(p,q)$ is the set of piecewise smooth curves from $p$ to $q$ and
$\ell_F(\gamma)$ is the Finslerian length of $\gamma:[a,b]\rightarrow \R$, that is
$$\ell_F(\gamma)=\int_a^b F(\dot\gamma(s))\df s.$$
Because of \eqref{triineq}, $\dist$  satisfies the triangle
inequality. As all the properties of a distance but symmetry are
fulfilled, the pair $(M,\dist)$ is referred sometimes as a {\em
generalized metric space} (see \cite[Section 6.2]{BCS01} or \cite{FHSpre}
for a detailed study). When the Finsler metric is non-reversible,
$\df$ is not symmetric, because the length of a curve $\gamma$ may
not coincide with the length of its reverse  curve
$\tilde{\gamma}(s)=\gamma(b+a-s)\in M$.
This also translates to Cauchy sequences, that is, we say that a
sequence $\{x_i\}_{i\in \N}$ in $M$  is {\em forward} (resp. {\em
backward}) Cauchy if for every $\varepsilon>0$ there exists $N\in
\N$ such that if $i,j>N$, then $\dist (x_i,x_j)<\varepsilon$
whenever $i\leq j$ (resp. $i\geq j$). Moreover, there are two
kinds of (open) balls, {\em forward} balls, that is,
$$B^+(p,r)=\{x\in M : \dist(p,x)<r\},$$
where $p\in M$ and $r\geq 0$, and {\em backward} balls,
$$B^-(p,r)=\{x\in M : \dist(x,p)<r\}.$$
As usual,  a bar  will denote closure,
 as in the closed balls  $\bar{B}^+(p,r)$ or $\bar{B}^-(p,r)$.
 The topologies generated by the forward and the backward balls
agree with the underlying manifold topology (see \cite[Section
6.2 C]{BCS01}).

\subsection{Geodesics.} There are several
ways to define geodesics of Finsler manifolds. We can use any of
the connections associated to a Finsler manifold, for example:
Chern, Cartan, Berwald or Hashiguchi connections (see
\cite{BCS01}). Another possibility is to define a geodesic as a
smooth critical curve of the length functional (see
\cite[Prop. 5.1.1]{BCS01}); nevertheless, as in the
Riemannian case, such a functional is invariant by
reparametrizations of the curves, and its critical points will be
called {\em pregeodesics} here\footnote{We  observe that we will
use a different convention from \cite{BCS01}, where pregeodesics
are called geodesics. So, our geodesics are the speed constant
geodesics in \cite{BCS01}.}. Geodesics affinely parametrized by
arc length are the critical points of the energy functional
\begin{equation}\label{E}
E_F(\gamma)=\int_a^b F^2(\dot\gamma(s))\df s,
\end{equation}
defined in the space of $H^1$-curves $\gamma:[a,b]\rightarrow M$
with fixed endpoints (see  for example \cite[Prop. 2.3]{cmm}).
Another consequence of non-reversibility is that geodesics are
non-reversible, that is,  the reverse curve of a geodesic may not
be a geodesic. This leads us to define two exponential maps at every
point $p\in M$; the first one will be taken as the natural
exponential, since it is analogous to the Riemannian exponential for
geodesics departing from $p$, say, ${\rm exp}_p(v)=\gamma_v(1)$,
where $\gamma_v$ is the unique geodesic, such that $\gamma_v(0)=p$
and $\dot\gamma_v(0)=v$. The second one, which we call {\it
reverse exponential map}, $\tilde{\rm exp}$, can be defined as the
exponential map associated to the  reverse Finsler metric
$\tilde{F}$  (see \cite[Chapter 6]{BCS01}).

\subsection{Geodesic completeness and Hopf-Rinow theorem.}
There are two types of geodesic completeness, {\it forward}  when
the domain of the geodesics can  be always extended to
$(a,+\infty)$, for some $a\in\R$, and {\it backward}, when it  can be
extended to $(-\infty, b)$ for some $b\in\R$. In these
circumstances, the classical Hopf-Rinow theorem splits into a
forward and a backward version (see \cite[Th. 6.6.1]{BCS01} or
\cite{Daz69}).
\begin{thm}\label{Hopf-Rinow}
Let $(M,F)$ be a Finsler manifold, then the following properties are equivalent:
\begin{itemize}
\item[(a)] The generalized metric space $(M,\dist)$ is forward (resp.
backward) complete. \item[(b)] The Finsler manifold $(M,F)$ is
forward (resp. backward) geodesically complete. \item[(c)] At
every point $p\in M$, ${\rm exp}_p$ (resp. $\tilde{\rm exp}_p$) is
defined on all of $T_pM$. \item[(d)] At some point $p\in M$, ${\rm
exp}_p$ (resp. $\tilde{\rm exp}_p$) is defined on all of $T_pM$.
\item[(e)] Heine-Borel property: every closed and forward (resp.
backward) bounded subset of $(M,\dist)$ is compact.
\end{itemize}
Moreover, if any of the above conditions holds, $(M,F)$ is {\em
convex}\footnote{Notice that there is no ``forward and backward''
convexity, as the former would be equivalent to the latter.},
i.e., every pair of points $p,q\in M$ can be joined by a
minimizing geodesic from $p$ to $q$.
\end{thm}
In order to overcome the lack of symmetry of the distance $\dist$ in \eqref{distance}, we can define the symmetrized distance as
\begin{equation}\label{symdist}
\ds(p,q)=\frac{1}{2}(\dist(p,q)+\dist(q,p)).
\end{equation}
It is easy to see that $\ds$  is a distance in the classical
sense, even though it is not constructed as a length metric (see
Appendix). We will denote its associated balls as $B_s(x,r)$,
for $x\in M$ and $r\geq 0$. We can wonder if a Hopf-Rinow theorem
holds also for $\ds$. As a first answer, we obtain the following
result.
\begin{prop}\label{dscomplete}
Let $(M,F)$ be a Finsler manifold. The following properties are
equivalent:
\begin{itemize}
\item[(a)] Heine-Borel property: the {\em symmetrized closed
balls} $\bar{B}_s(x,r)$ are compact for all $x\in M$ and  $r>0$.
\item[(b)] $\bar{B}^+(x,r)\cap \bar{B}^-(x,r) $ is compact for any
$x\in M$ and  $r>0$. \item[(c)]  $\bar{B}^+(x,r_1)\cap
\bar{B}^-(y,r_2) $ is  compact for any $x,y\in M$ and $r_1,r_2>0$.
\end{itemize}
Moreover, if any of the above conditions hold, then the metric
space  $(M,\ds)$ is complete.
\end{prop}
\begin{proof} (a) $\Rightarrow$ (b). It follows from the obvious
inclusion $\bar{B}^+(x,r)\!\!\cap \bar{B}^-(x,r)\subset
\bar{B}_s(x,r)$, as the latter ball is compact.

(b) $\Rightarrow$ (c). As $\bar{B}^-(y,r_2)\subset
\bar{B}^-(x,r_2+\dist(y,x))$, the result follows from
$\bar{B}^+(x,r_1)\cap \bar{B}^-(y,r_2)\subset \bar{B}^+(x,r_3)\cap
\bar{B}^-(x,r_3)$, where $r_3=\max\{r_1,r_2+\dist(y,x)\}$ .

(c) $\Rightarrow$ (a). Straightforward from  $\bar{B}_s(x,r)\subset
\bar{B}^+(x,2r)\cap \bar{B}^-(x,2r)$.

Finally, $(\rm a)$ implies that the metric $\ds$ is complete
because a  $\ds$-Cauchy sequence is always contained in a ball
$\bar{B}_s(x,r)$ for some $x\in M$ and $r>0$.
\end{proof}
Th. \ref{Hopf-Rinow} cannot be claimed to prove convexity
under the Heine-Borel property in the last proposition, as hypotheses
$(\rm a)-(\rm c)$ above are weaker than those in the theorem.
 However, we will see in Section
\ref{randersfermat} that these hypotheses do imply convexity  (see
Th. \ref{tgeodesicconnection}). Nevertheless, there is no relation
between the hypotheses and  completeness. In fact, Example
\ref{minkexe} shows a Randers  metric that satisfies $(\rm
a)$-$(\rm c)$ in Prop. \ref{dscomplete}, but with forward and
backward incomplete geodesics. Moreover, in the following example
we exhibit a Randers metric with complete symmetrized distance
that does not satisfy the equivalent conditions $(\rm a)$ to $(\rm
c)$ of Prop. \ref{dscomplete} (in fact, the manifold is
$\ds$-bounded but non-compact), see also Subsection \ref{s5.3}.

\begin{exe}\label{counterhopfrinow}
Consider $\R^2$ and there two  smooth  bump functions $\mu_+,
\mu_-$  such that $0\leq\mu_\pm\leq 1$ and satisfying :
$$
\mu_\pm(x,y)\equiv \mu_\pm(x)=\left\{
\begin{array}{llll}
0& & \mbox{if} & |x\mp 3|\geq 2 \\
1& & \mbox{if} & |x\mp 3|\leq 1
\end{array}
\right.
$$
and the metric $h=\dist x^2+\dist y^2.$
Now, consider the $1$-form $\omega$:
$$
\omega=(\mu_+(x)-\mu_-(x))\frac{y^2}{1+y^2}\dist y
$$
and the corresponding Randers metric: $F(v)= \sqrt{h(v,v)}+
\omega(v)$ for $ v\in T\R^2$, with associated distance $\dist$
and symmetrized distance $\ds$.
 Finally, consider the
strip $-6\leq x \leq 6$, construct a quotient $M$ by identifying
each two $(-6,y), (6,y)$, and regard $F$ as a Randers metric on
$M$.

Easily,  the lines $ s\mapsto (3,s)$ and $ s\mapsto (-3,-s)$
have infinite length for $F$, whereas the curves $s\mapsto (3,-s)$
and $ s\mapsto (-3,s)$
 have finite length equal to $\pi$. Thus, the
distance $\dist$, and therefore the symmetrized one $\ds$, are
finite (say, obviously bounded by  $12+\pi$). As a
consequence, neither $\dist$ nor $\ds$ can satisfy the property of
Heine-Borel ($M$ is non-compact but included in a ball), and
$\dist$ is incomplete. Nevertheless, $\ds$ is complete, as no
non-converging sequence $\{p_n\}_n$ can be Cauchy; in fact, either
the limsup of $\{\dist(p_m,p_{m+k})\}_k$ or the one of
$\{\dist(p_{m+k},p_m)\}_k$ is bounded away from zero.
\end{exe}
\section{Causality of spacetimes}\label{causality}
\subsection{Lorentzian manifolds and spacetimes.}
Our notation and conventions on Causality will be standard, as in
\cite{BeErEa96, HawEll73, MinSan06, oneill83, SaWu77}. So, for a
Lorentzian manifold
 $(\GG, \g)$,  ${\dim M}\geq 2$,
the metric $g$ on $\GG$ has index  $(-,+,\dots,+)$, and a tangent
vector $v\in T\GG$, is {\it timelike} when
$g(v,v)<0$, {\it spacelike} when $g(v,v)>0$, {\it lightlike} if
$g(v,v)=0$
 but $v\not=0$ and {\it causal} if it is timelike or lightlike; following \cite{MinSan06}
  vector 0 will be regarded as non-spacelike and non-causal --even though this is not by any means the unique convention  in the literature.
  At every point $p\in \GG$  the {\em causal cone} is the subset of causal vectors in $T_p\GG$, which
  has exactly two connected components.
  A {\em time-orientation} is a smooth
  choice of a causal cone at every point, which will be called the {\em future} causal cone --in opposition to the non-chosen one or {\em past} causal cone.
    A {\it spacetime} is a connected $C^\infty$  Lorentzian manifold $(\GG,  \g)$
    endowed  with a time-orientation  (which is not written explicitly in the notation). The latter can be determined by a timelike vector
    field $T$ which defines the future orientation and, so, a causal vector $v\in T\GG$ is future-pointing (resp. past-pointing) if
    $g(v,T)<0$ (resp. $g(v,T)>0$).
    A piecewise smooth curve $\gamma:[a,b]\rightarrow \GG$ will be said timelike (and analogously spacelike, lightlike, causal, or future/past-pointing)
    if so is  its velocity $\dot\gamma(s)$  at every $s\in[a,b]$.
Spacetimes are used in {\it General Relativity} as models of
(regions of the ) Universe. The points of $M$ are also called
{\em events} (they represent all possible ``here-now'') and
massive (resp. massless) particles are
 described by  future-pointing timelike (resp. lightlike)
 curves.

Causality studies the properties associated to the causal cones,
as, for example, if two  events can be
connected by means of a causal curve. As two Lorentzian metrics
on the same manifold are (pointwise) conformal iff they have equal
causal cones,
Causality is essentially the same thing as
conformal geometry in Lorentzian Geometry (even though usually the
former refers to the global viewpoint of the latter). Given two
events $p$ and $q$ in a spacetime, we say that they are
chronologically related, and write $p\ll q$ (resp. strictly
causally related $p< q$)  if there exists a future-pointing
timelike (resp.
 causal) curve $\gamma$ from $p$ to $q$; $p$ is causally related to $q$ if either $p<q$ or $p=q$, denoted $p\leq q$.
Relations such as $p\leq q \ll r \Rightarrow p \ll r$ are well-known.
The {\it chronological future}  (resp. {\it causal future}) of
$p\in \GG$ is defined as $I^+(p)=\{q\in \GG : p\ll q\}$ (resp.
$J^+(p)=\{q\in \GG : p\leq q\}$). Analogous notions appear
reversing the word ``future" by ``past" and, so, one writes
$I^-(p), J^-(p)$.

\subsection{Causal properties of spacetimes.}

The {\em causal ladder} groups spacetimes in the following
families,  ordered by strictly increasingly better causal
properties:
\begin{multline*}
\text{chronological}\Leftarrow
\text{causal}\Leftarrow\text{distinguishing}\Leftarrow\text{strongly
causal}\\
\Leftarrow\text{stably causal}\Leftarrow\text{causally
continuous}\\\Leftarrow\text{causally
simple}\Leftarrow\text{globally hyperbolic}
\end{multline*}
In the following we will give a brief account of these spacetime
classes; for further information see
\cite{BeErEa96,HawEll73,SaWu77}, or the  survey \cite{MinSan06}.

A spacetime is {\it chronological} (resp. {\em causal}) if
$p\notin I^+(p)$ (resp. $p\notin J^+(p)$) for every $p\in \GG$;
this comprises the nonexistence of timelike or causal closed
curves. A spacetime is future (resp. past) {\it distinguishing} if
$I^+(p)=I^+(q)$ (resp. $I^-(p)=I^-(q)$) implies $p=q$, and {\it
distinguishing} if it is both future and past distinguishing. It
is easy to prove that distinguishing spacetimes are causal.
Intuitively, {\em strong causality} means the in existence of
almost closed timelike curves, and this is equivalent to obtaining
a basis of the manifold topology with the subsets  of the type
$I^+(p)\cap I^-(q)$. A spacetime is {\it stably causal} if it is
causal and it remains causal when we open slightly the light
cones. This is equivalent (see \cite{BeSan05} or \cite{San05}) to
the existence of a temporal function on $(\GG, \g)$, that is, a
smooth function $t:\GG\rightarrow\R$ with a past-pointing timelike
gradient (thus, $t$ is also  a {\em time function}, i.e.,
continuous and strictly increasing on every future-pointing causal
curve). A spacetime is said  causally continuous when the maps
$I^\pm:\GG\rightarrow {\mathcal P}(\GG)$ are one to one (i.e., the
spacetime is distinguishing) and continuous (here ${\mathcal
P}(\GG)$ is the set of parts of $\GG$ endowed with the topology
which admits as  a basis  the collection
$\{\mathcal{O}_K\}_{K\subset \GG}$, where each open
$\mathcal O_K$ contains all the subsets of $\GG$ which do not
intersect the compact set $K\subset \GG $, see \cite[Defn. 3.59,
Prop. 3.38]{MinSan06}). A spacetime is {\it causally simple} when
it is causal and all causal  futures and pasts $J^\pm(p)$ are
closed for every $p\in \GG$ \cite{BeSan07}.
Finally, a spacetime is {\it globally hyperbolic} when it admits a
Cauchy hypersurface, that is, a subset $S$ which meets exactly
once every inextensible timelike curve --which can be chosen as  a
smooth spacelike hypersurface, necessarily crossed once by
inextensible causal curves. This is equivalent to be causal with
$J^+(p)\cap J^-(q)$ compact for every pair $p,q\in \GG$ (see
\cite{BeSan03, BeSan05, BeSan06, BeSan07,Ge70} or \cite[Sections
3.11.2, 3.11.3]{MinSan06} for details).

\section{Fermat metrics applied to stationary spacetimes }\label{stationaryfermat}

A spacetime $(M,g)$ is called {\it stationary} if it admits a {\em
stationary vector field}, i.e., a timelike Killing vector field.
Let $\R\times S$ be a product manifold with natural projection on
the first factor  $ t:\R\times S\rightarrow \R$, called {\em standard time}, and on
the second factor $\pi: \R\times S\rightarrow S$. $(\R\times S ,
g)$ is a {\em standard stationary} spacetime if there exists a
Riemannian metric $g_0$ on $S$, a positive function
$\beta:S\rightarrow \R$ and a smooth  $1$-form
$\omega$ on $S$, such that:
\begin{equation}\label{standard}
 \g =-(\beta\circ\pi)  dt^2+\pi^*\omega\otimes  d t + d t \otimes \pi^*\omega+\pi^*g_0,
\end{equation}
where $\pi^*\omega$ and $\pi^* g_0$ are the pullback of $\omega$ and $g_0$ on $\R\times S$ through $\pi$. The (future) time-orientation
 is determined by the {\em standard stationary vector field} $K$
 obtained as the lift  to $\R\times S$ through
$ t$
of the natural vector field on  $\R$. 

A stationary spacetime with a prescribed stationary vector field
$X$ admits an isometry with a standard stationary spacetime whose
differential maps $X$ in $K$, if and only if the flow of $X$ is
complete and $M$ admits a  spacelike section i.e., a
spacelike hypersurface $S$ which is crossed exactly once by any
inextensible integral curve of $X$. In this case, an isometry with
(\ref{standard}) is obtained by moving $S$ by means of the flow of
$X$, and all the possible isometries correspond with all the
possible spacelike sections. Remarkably, the existence of such a
section is determined by the level where the spacetime is positioned  in the causal
ladder. Concretely, a stationary spacetime is isometric to a
standard stationary one if and only if  it is distinguishing and
admits a complete stationary vector field (see \cite{JaSan08} for
this result and other details).


\subsection{Fermat principle.}\label{fermatprinciple}
According to the relativistic Fermat principle, any lightlike
pregeodesics is a critical point of the arrival time function
corresponding to an {\em observer} (defined as a future-pointing
timelike curve, up to re-parametrization), see
\cite{Kovner90,Per90}. In a standard stationary spacetime
$(\R\times S, \g)$, this implies that future-pointing  lightlike
geodesics project onto pregeodesics of a Finsler metric in $S$.
More precisely, when a vertical line $\R\ni s\mapsto
(s,x_1)\in\R\times S$ is regarded as an observer,  the arrival
time AT($\gamma$) of a (future-pointing) lightlike curve
$\gamma:[a,b]\rightarrow \R\times S$, $\gamma=(t_\gamma,x_\gamma)$
(i.e. $\gamma(s)=(t_\gamma(s),x_\gamma(s)), s\in [a,b]$), joining
a point $(t_\gamma(a),x_\gamma(a))$ with the vertical line, can be
expressed as ${\rm AT}(\gamma)=t_\gamma (b)$, where
$t_\gamma(b)-t_\gamma(a)$ is equal to:
\begin{equation}\label{arrivaltime}
\int_a^b \left(\frac{1}{\beta(x_\gamma)}\omega(\dot
x_\gamma)+\sqrt{\frac{1}{\beta(x_\gamma)}g_0(\dot x_\gamma,\dot x_\gamma)+\frac{1}{\beta(
x_\gamma)^2} \omega(\dot x_\gamma)^2}\right)\; \df s.
\end{equation}
As a result, future-pointing lightlike geodesics  project  onto
the pregeodesics of the non-reversible Finsler metric of Randers
type  on $S$
\begin{equation}\label{fermat}
F(v)=\frac{1}{\beta}\omega(v)+\sqrt{\frac{1}{\beta}g_0(v,v)+\frac{1}{\beta^2}\omega(v)^2}, \quad  \forall v\in TS
\end{equation}
which we call  {\em Fermat metric}. Prescribing $x_\gamma$, the
$t_\gamma$-component can be recovered  as $t_\gamma(s)=
t_\gamma(a)+ \int_a^sF(\dot x_\gamma(\mu))\df \mu$.
By a similar reasoning, past-pointing lightlike geodesics project
onto pregeodesics of the reverse Finsler metric
\begin{equation*}
\tilde{F}(v)=-\frac{1}{\beta}\omega(v)+\sqrt{\frac{1}{\beta}g_0(v,v)+\frac{1}{\beta^2}\omega(v)^2}.
\end{equation*}
Such Fermat metrics were introduced in \cite{cmm} to obtain some
multiplicity results for lightlike geodesics and timelike
geodesics with fixed proper time from an event to a vertical line
in globally hyperbolic standard stationary spacetimes as a
consequence of multiplicity results for Finsler metrics.

Prescribing a piecewise smooth  future-pointing lightlike curve
$\gamma:[a,b]\rightarrow\R\times S, \gamma=(t_\gamma,x_\gamma)$,
its projection on $S$ is a curve with Fermat length
$\ell_F(x_\gamma)=t_\gamma(b)-t_\gamma(a)$. Conversely, a
piecewise smooth curve $x_\gamma:[a,b]\rightarrow S$ can be lifted
to a  future-pointing lightlike curve $[a,b]\ni
s\rightarrow\gamma(s)=(t_\gamma(s),x_\gamma(s))\in \R\times S$ by
choosing $t_\gamma(s)=\int_a^sF(\dot
x_\gamma(\mu))\df \mu$ and therefore,
$t_\gamma(b)-t_\gamma(a)=\ell_F(x_\gamma)$. Easily then:
\begin{prop}\label{ppp}
Let $z_0=(t_0,x_0)$, $L_{x_1}=\{(t,x_1): t\in \R\}$ be,
respectively, a point and a line in a standard stationary
spacetime. Then $z_0$ can be joined with $L_{x_1}$ by means of a
future-pointing (resp. past-pointing) lightlike pregeodesic
$t\mapsto\gamma(t)=(t,x_\gamma(t))$ (resp.
$t\mapsto\gamma(t)=(2t_0-t,x_\gamma(t))$ starting at $z_0$, if and
only if $x_\gamma$ is a unit  speed geodesic of the Fermat metric
$F$ (resp. $\tilde F$) which joins $x_0$ with $x_1$. In this case,
the interval of time $t_1-t_0$ (resp. $t_0-t_1$) such that
$\gamma(t_1)\in L_{x_1}$ is equal to the length of the curve
$x_\gamma(t), t\in [t_0,t_1]$ (resp. $t\in [t_0, 2t_0-t_1]$)
computed with $F$ (resp. $\tilde F$).
\end{prop}
 In fact, a (future-pointing) curve parametrized with the standard time $t\mapsto \gamma(t)=(t,x_\gamma(t))$ is lightlike iff $x_\gamma$
 is parametrized with Fermat speed  $F(\dot x_\gamma)=1$  and, among these curves,
 Fermat's principle states that the critical curves of the arrival time \eqref{arrivaltime}
 are lightlike pregeodesics.
 As the arrival time coincides (up to  an initial additive constant)
 with the Fermat length of $x_\gamma$, it follows that $x_\gamma$ must be a geodesic for the Fermat metric   (see also
\cite[Th. 4.1]{cmm}).
\subsection{Causality via Fermat metrics.} The essential point
 about Fermat metrics is that they contain all the causal information of a standard stationary spacetime.
Let $\gamma=(t_\gamma,x_\gamma):[0,1]\rightarrow \R\times S$, be a
future-pointing differentiable causal curve, then
\begin{equation*}
 g_0(\dot x_\gamma,\dot x_\gamma)+2\omega(\dot x_\gamma)\dot t_\gamma-\beta (x_\gamma) \dot t_\gamma^2\leq 0.
\end{equation*}
Analyzing the quadratic equation in $\dot t_\gamma$ and using that
$\beta>0$ we deduce that either
\begin{align}
& \dot t_\gamma\geq \frac{1}{\beta (x_\gamma)}\omega(\dot x_\gamma)+\sqrt{\frac{1}{\beta(x_\gamma)}\omega(\dot x_\gamma,\dot x_\gamma)+\frac{1}{\beta(x_\gamma)^2}\omega(\dot x_\gamma)^2}=F(\dot x_\gamma)\nonumber
\intertext{or}
& \dot t_\gamma\leq \frac{1}{\beta(x_\gamma)}\omega(\dot
x_\gamma)-\sqrt{\frac{1}{\beta(x_\gamma)}g_0(\dot x_\gamma,\dot x_\gamma)+\frac{1}{\beta(
x_\gamma)^2}\omega(\dot x_\gamma)^2}=-
\tilde{F}(\dot x_\gamma).\label{segundat}
\end{align}
As  $\gamma$ has been chosen to be future-pointing,
 we have  $ \g(K,(\dot t_\gamma,\dot x_\gamma))=\omega(\dot
x_\gamma)-\beta(x_\gamma) \dot
t_\gamma<0$. This implies that \eqref{segundat} cannot be
satisfied and then $\dot t_\gamma\geq F(\dot
x_\gamma)\geq 0$. Moreover, $\dot t_\gamma>0$ because, otherwise,
$\dot t_\gamma(s)=0=\dot x_\gamma(s)$ at some $s$ and $\gamma$
would not be causal there (in particular, this proves that the
standard time $t:\R\times S\rightarrow \R$ is a temporal function
and the spacetime is stably causal). In \cite[Th. 4.2]{cmm} it is
proven that
 if the Fermat metric is forward
or backward complete, then the spacetime is globally hyperbolic
and if the slice $S_0=\{0\}\times S$ is a Cauchy hypersurface then
the Fermat metric is forward and backward complete. In the
following we will give a complete characterization of the causal
properties in terms of the Fermat metric. As a first step,
 let us obtain a description of the chronological past and
future.
 From now on, the balls $B^+(x,r)$ and $B^-(x,r)$ will correspond
to the forward and the backward balls in the generalized metric space
$(S,F)$, $F$ defined in \eqref{fermat}.
\begin{prop}\label{bolas}
Let $(\R\times S,  \g)$ be a standard stationary spacetime as in
\eqref{standard}. Then
 $$I^+(t_0,x_0)= \bigcup_{s> 0}\{t_0+s\}\times B^+(x_0,s), $$
 $$I^-(t_0,x_0)= \bigcup_{s< 0}\{t_0-s\}\times B^-(x_0,s).$$
\end{prop}
\begin{proof}
(We will consider just the first equality.) Let $(t_1,x_1)\in
I^+(t_0,x_0)$. As
$I^+(t_0,x_0)$
is open, one finds easily a
lightlike piecewise geodesic $\gamma=(t_\gamma,x_\gamma)$ joining
$(t_0,x_0)$ and
$(t_1-\varepsilon,x_1)$
for some small $\varepsilon>0$ such that
$(t_1-\varepsilon,x_1)$
remains in
$I^+(t_0,x_0)$
 (see for example \cite[Prop. 2]{FloSan08}). Then $x_\gamma$
is a curve joining $x_0$ and $x_1$ with Fermat length equal to
$t_1-t_0-\varepsilon$, so that $x_1\in B^+(x_0,t_1-t_0)$ and
$(t_1,x_1)\in \{t_0+(t_1-t_0)\}\times B^+(x_0,t_1-t_0)$.

Conversely, let $x_1\in B^+(x_0,s)$, and take an arc length
parametrized curve $x_\gamma:[0,b]\rightarrow  S$ joining $x_0$
and $x_1$ with Fermat length $b<s$. The future-pointing lightlike
curve $\gamma(r)= (t_0+r, x_\gamma(r))$ for $r\in[0,b]$  yields
$(t_0,x_0)\leq (t_0+b,x_1)$. As, trivially,   $(t_0+b,x_1)\ll
(t_0+s,x_1)$ we
 have
 $(t_0+s,x_1)\in  I^+(t_0,x_0)$.
\end{proof}
\begin{thm}\label{causalstructure}
Let $(\R\times S,  \g)$ be a standard stationary spacetime. Then
$(\R\times S,  \g)$ is causally continuous and
\begin{itemize}
 \item[(a)]
 the following assertions become equivalent:
\begin{itemize}
\item[(i)]
$(\R\times S,  \g)$ is causally simple, \item[(ii)] $J^+(p)$ is closed
for all $p$, \item[(iii)] $J^-(p)$ is closed for all $p$ and \item[(iv)] the
associated Finsler manifold $(S,F)$ is convex,
\end{itemize}
\item[(b)] it is globally hyperbolic  if and only
if the symmetrized closed balls $\bar B_s(x,r)$ are compact for
every $x\in S$ and $r>0$.
\end{itemize}
\end{thm}
\begin{proof} The first assertion is known (see \cite{JaSan08})
but we can give a simple proof in terms of Fermat metrics. It is
enough to prove that $(\R\times S, \g)$ is future and past
reflecting (see for example \cite[Defn. 3.59, Lemma
3.46]{MinSan06} or \cite[Th. 3.25, Prop. 3.2]{BeErEa96}), and we
will focus on the latter,  that is, $I^+(p)\supset I^+(q)$ implies
$I^-(p)\subset I^-(q)$.  Let $p=(t_0,x_0)$ and $q=(t_1,x_1)$. Then
$I^+(p)\supset I^+(q)$ implies that ${\dist}(x_0,x_1)\leq t_1-t_0$
 by the first equality in Prop. \ref{bolas}.  Therefore, using the second equality,
$I^-(p)\subset I^-(q)$.

For the remainder, put $p=(t_0,x_0)$ and recall that,  by using
Prop. \ref{bolas} and \cite[Lemma 14.6]{oneill83} resp.:
\begin{equation}\label{eaux} \overline{I}^+(t_0,x_0)=
\bigcup_{s\geq 0}\{t_0+s\}\times \bar B^+(x_0,s) \quad \hbox{and}
\quad J^+(p) \subset \overline{ I}^+(p)= \overline{ J}^+(p).
\end{equation} 
Notice also that the condition of {\em causality} in
the definitions of causal simplicity and global hyperbolicity are
automatically satisfied.

For the proof of (a), it is enough to check that (ii)
and (iv) are equivalent.

Implication (ii) $\Rightarrow$ (iv).
Take any pair of points $x_0, x_1$ in $S$. From (\ref{eaux}),
$$({\rm d}(x_0,x_1),x_1)\in J^+(0,x_0)\setminus
I^+(0,x_0),$$  and  there exists a future-pointing lightlike
geodesic joining the points $(0,x_0)$ and $(\dist(x_0,x_1),x_1)$ (see for
example \cite[Prop. 10.46]{oneill83}). Clearly, the projection in
$S$ of this geodesic is the
required minimizing Fermat geodesic.

Implication  (ii) $\Leftarrow$ (iv). We have to prove that the
inclusion in (\ref{eaux}) is an equality. Any $(t_1,x_1)\in
\partial J^+(p)$ satisfies $t_1=t_0+\dist(x_0,x_1)$. So, take the
minimal Fermat geodesic $x$ starting at $x_0$ and ending at $x_1$.
The associated lightlike geodesic starting at $p$ in $(\R\times S,
\g)$ connects $p$ and $(t_1,x_1)$ as required.

Proof of (b). ($\Rightarrow$)
Consider the points $(r,x)$ and $(-r,x)$. Global hyperbolicity
implies  that $J^\pm(p)$ is closed; thus, Prop. \ref{bolas}  and
(\ref{eaux}) yield:
$$
 \{0\}\times
\left(\bar{B}^+(x,r)\cap \bar{B}^-(x,r)\right)= (\{0\}\times S)
\cap J^+(-r,x)\cap J^-(r,x),$$ and the right-hand side is compact,
also by global hyperbolicity. Then by Prop.~\ref{dscomplete} we conclude.

($\Leftarrow$)
Given two points $(t_0,x_0)$ and $(t_1,x_1)$ in $\R\times S$,
\begin{equation}\label{closurecompact}
J^+(t_0,x_0)\cap J^-(t_1,x_1)\subset
\bigcup_{s\in[0,t_1-t_0]}\{t_0+s\}\times\left(\bar{B}^+(x_0,s)\cap\bar{B}^-(x_1,t_1-t_0-s)\right).
\end{equation}
Moreover, the subset in the right-hand side is compact. Indeed,
any sequence $\{(s_k,y_k)\}_k$ in it has $\{s_k\}\subset
[0,t_1-t_0]$ and, thus, $\{s_k\}\rightarrow \bar{s}$, up to a
subsequence. Moreover, $\{y_k\}_k\subset
\bar{B}^+(x_0,t_1-t_0)\cap\bar{B}^-(x_1,t_1-t_0)$, which is
compact (see part $(\rm c)$ of Prop. \ref{dscomplete}). Thus,
again, $\{y_k\}\rightarrow \bar{y}$, up to a subsequence, and by
the continuity of the distance $(\bar{s},\bar{y})\in \{t_0+
\bar{s}\}\times
\bar{B}^+(x_0,\bar{s})\cap\bar{B}^-(x_1,t_1-t_0-\bar{s})$, which
concludes the compactness. By Eq. \eqref{closurecompact}, the
closure of $J^+(t_0,x_0)\cap J^-(t_1,x_1)$ is compact. As the
spacetime is strongly causal, this implies that $J^+(t_0,x_0)\cap
J^-(t_1,x_1)$ is compact  (see \cite[Lemma 4.29]{BeErEa96}), and
global hyperbolicity follows.
\end{proof}

Th. \ref{causalstructure} allows us to determine easily examples
of standard stationary spacetimes which are not causally simple,
even in the static case  ($\omega =0$), extending \cite[Rem.
3.2]{SaNonLin05}. Notice also that we have characterized global
hyperbolicity, which is a property intrinsic to the spacetime,
independently of how it is written as a standard stationary one.
Nevertheless, the fact that $S$  is a Cauchy hypersurface (more
precisely, a slice $\{t_0\}\times S$, and trivially then any slice,
is Cauchy) will depend on the concrete choice, and it is
characterized next.

\begin{thm}\label{th2}
Let $(\R\times S,  \g)$ be a standard stationary spacetime.  A
slice $S_{t_0}=\{t_0\}\times S,\ t_0\in\R$, is a Cauchy hypersurface if and
only
if the Fermat metric $F$ on $S$
is forward and backward complete.
\end{thm}
\begin{proof}
As the slice is spacelike and acausal, it  is Cauchy iff any future-pointing
inextensible null pregeodesic $\gamma\colon(a,b)\rightarrow
\R\times S$ meets  $S_{t_0}$  once (see \cite[Prop.
5.14]{Penros72} or \cite[Cor. 14.54]{oneill83}) and, in this case,
$\gamma$ crosses all the slices. As $\gamma$ can be parametrized
with the standard time, $\gamma(t)=(t,x_\gamma(t))$, this curve
will cross all the slices iff the inextensible domain $(a,b)$ of
its $x_\gamma$-component is equal to $\R$. But the possible
$x_\gamma$-components are all the (unit speed) Fermat geodesics,
so, their domains are equal to $\R$ if and only if $(S,F)$ is
forward and backward complete.
\end{proof}

\begin{rem} \label{r1} From the proof of Th. \ref{th2}, a
more precise result follows:

{\em $(S,F)$ is forward (resp. backward)  complete  iff any
future-pointing (resp. past-pointing) inextensible lightlike
geodesic --and, then, also any timelike curve-- starting at
$S_{0}=\{0\}\times S$ crosses all the slices $S_t$ for $t>0$ (resp.
$t<0$).}

Informally, this means that forward/backward completeness is
equivalent to the property that the slices behave as  Cauchy
hypersurfaces for  future/past-pointing causal curves. As a
consequence, one can construct a Fermat metric, with compact
symmetrized closed balls, which is forward and (or) backward
incomplete, by taking a globally hyperbolic stationary spacetime
and splitting it as a standard one with respect to a non-Cauchy
spacelike hypersurface $S$. The following example illustrates this
situation.
\end{rem}

\begin{exe}\label{minkexe}
Consider Lorentz-Minkowski spacetime $\mathbb{L}^2$, i.e $\R^2$ endowed with the metric
$g=-\dist t^2 + \dist x^2$  and the spacelike section given by the
curve
\[\alpha(\theta)=\left\{
\begin{array}{ll}
(-\cosh \theta +1,\sinh \theta)
 & \text{if }\quad \theta\in(-\infty,0],\\
(\cosh \theta-1,\sinh \theta) & \text{if }\quad
\theta\in[0,+\infty).
\end{array}\right.
\]
As emphasized at the beginning of  the present section, we can
express the Minkowski metric as a standard stationary spacetime by
using the ``spacelike hypersurface'' $\alpha$ (one can easily smooth $\alpha$ in a neighborhood of $0$, this will not affect the
discussion below).  Putting $
v_\theta=\alpha'(\theta)$, a new standard stationary splitting of
$\mathbb{L}^2$ is determined by $\beta\equiv 1$, $g_0(\mu
v_\theta,\mu v_\theta )=\mu^2$ and
\[
\omega(\mu v_\theta)=g(\partial_t,\mu v_\theta)/\beta =
\left\{
\begin{array}{lll}
\mu \sinh \theta &\text{if} &\theta\leq 0,\\
-\mu \sinh \theta& \text{if}& \theta\geq 0.
\end{array}\right.
\]
The associated Fermat metric is
\[F(\mu v_\theta)=\left\{
\begin{array}{ll}
\sqrt{\mu^2(1+\sinh^2\theta)}+\mu\sinh \theta
 & \text{if }\quad \theta\in(-\infty,0],\\
\sqrt{\mu^2(1+\sinh^2\theta)}-\mu\sinh \theta
& \text{if }\quad \theta\in[0,+\infty).
\end{array}\right.
\]
The length of $\R$ with this metric is
\[\int_{-\infty}^{+\infty}F(v_\theta)\dist\theta =
\int_{-\infty}^0(\cosh\theta+\sinh\theta)\dist\theta+\int_0^{+\infty}(\cosh\theta-\sinh\theta)\dist\theta=2\]
and, thus $(\mathbb{R},F)$ is neither forward nor backward
complete, even though its symmetrized closed balls are compact.
Obvious modifications in the branches of $\alpha$ yield  only
forward and backward completeness, as pointed out in Rem.
\ref{r1}.
\end{exe}

\subsection{Cauchy developments.}\label{s4.3}
As the final application in this section, we will construct Cauchy developments in terms of the Fermat metric. A subset $A$ of a
spacetime $M$ is {\it achronal} if no  $x,y\in A$ satisfies $x\ll
y$; in this case,  the {\it future (resp. past) Cauchy
development} of $A$, denoted by $D^+(A)$ (resp. $D^-(A)$), is the
subset of points $p\in M$  such that every  past- (resp. future)-
inextensible causal curve through $p$ meets $A$. The union
$D(A)=D^+(A)\cup D^-(A)$ is the {\it Cauchy development} of $A$.
The {\it future $H^+(A)$ (resp. past $H^-(A)$) Cauchy horizon} is
defined as
\[H^\pm(A)=\{p\in \bar{D}^{\pm}(A):I^\pm(p) \cap D^\pm(A)=\emptyset\}.\]
Intuitively, $D(A)$ is the region of $M$ a priori predictable from
data in $A$,  and its  {\em horizon} $H(A)=H^+(A)\cup H^-(A)$, the
boundary of this region.

\begin{prop}\label{p4.5}
Let $(\R\times S,  \g)$ be a standard stationary spacetime as in
\eqref{standard}, with $S_0=\{0\}\times S$ a Cauchy hypersurface,  $A\subset
S$, and $A_{t_0}= \{t_0\}\times A$  the corresponding (necessarily
achronal) subset of $\{t_0\}\times S$. Then
\begin{equation}\label{cauchydevelop}
D^+(A_{t_0})=\{(t,y):\df (x,y)>t-t_0\,\,\text{for every $x\notin
A$ and $t\geq t_0$}\},
\end{equation}
\begin{equation}\label{cauchydevelop-}
D^-(A_{t_0})=\{(t,y):\df (y,x)>t_0-t\,\,\text{for every $x\notin
A$ and $t\leq t_0$}\},\end{equation} where $\df$ is the distance
in $S$ associated to the Fermat metric.

Moreover, the Cauchy horizons can be described as
\begin{equation}\label{horizon}
H^+(A_{t_0})=\{(t,y):\inf_{x\notin A}\df (x,y)=t-t_0
\}
\end{equation}
\begin{equation}\label{horizon-}
\quad H^-(A_{t_0})=\{(t,y):\inf_{x\notin A}\df (y,x)=t_0-t
\}.
\end{equation}
\end{prop}
\begin{proof}
Clearly, $D^+(A_{t_0})$ is contained in the semi-space
$t\geq t_0$.
Given a point $(t,y)\in \R\times S,\ t\geq t_0$,  every
past-inextensible causal curve meets the Cauchy hypersurface
$\{t_0\}\times S$ at some point. From the definition, $(t,y)
\not\in D^+(A_{t_0})$ iff there exists $x\in A^c$  (where $A^c$ is
the complementary subset of $A$ in $S$) such that $(t,y)\in
J^+(t_0,x)$. As $J^+(t_0,x)$ is closed, this is equivalent to
$d(x,y)\leq t-t_0$ (recall (\ref{eaux})), and the conclusion on
$D^+(A_{t_0})$ follows.

The characterization of the Cauchy horizon is obtained by taking into
account that  $\bar{D}^+(A_{t_0})=\{(t,y):\inf_{x\notin A}\df
(x,y)\geq t-t_0;\ t\geq t_0\}$ and using
 the property
 $(s,x)\in I^+(t,y)$ iff $\df
(y,x)<s-t$, which follows from  Prop. \ref{bolas}. Assume that
$\inf_{x\notin A}\dist(x,y)=t-t_0$, and thus $(t,y)\in \bar
D^+(A_{t_0})$. Then by the Finslerian Hopf-Rinow theorem this
infimum
is attained at some $\bar{x}\in \overline{A^c}$,
 i.e.  $\dist (\bar{x},y)= t-t_0$.
Let $(s_0,x_0)\in I^+(t,y)$, then $\dist (\bar{x},x_0)\leq
\dist(\bar x,y)+\dist(y,x_0)< t-t_0+s_0-t = s_0-t_0$. Moreover, as
$\bar{x}\in \overline{A^c}$ there exists $\tilde{x}\in A^c$ such
that $\dist (\tilde{x},x_0)<s_0-t_0$. Therefore,
\eqref{cauchydevelop} implies  that $(s_0,x_0)\notin
D^+(A_{t_0})$, and the inclusion $\supset$ in \eqref{horizon} is
proved. The other inclusion follows easily. Indeed, if there is
$(t,y)\in H^+(A_{t_0})$ such that $\inf_{x\notin
A}\dist(x,y)>t-t_0$, then we can choose $\varepsilon>0$ small
enough such that $(t+\varepsilon,y)\in D^+(A_{t_0})\cap I(t,y)$.
The pasts are obtained analogously.
\end{proof}

\begin{rem} Such a
 result can be extended in some different directions:

 (A) If $\R\times S$ is  globally hyperbolic, then any acausal compact
spacelike submanifold $A$ with boundary can be extended to a
spacelike Cauchy hypersurface $S_A$ \cite{BeSan06}.  Thus, $D(A)$
can be computed in terms of the Fermat metric associated to the
standard stationary  splitting for $S_A$.

(B) Even if $\R\times S$ is not globally hyperbolic (or $S$ is not
Cauchy) Cauchy developments could be studied by using the  Cauchy
boundary associated to the Finslerian metric (see \cite{FHSpre},
for properties of this boundary).
\end{rem}

Next, the results in \cite{LiNiren} on the regularity of the
Finslerian distance function  from the boundary will be used to
obtain some extensions of the results on differentiability  of
horizons in \cite{BK98} for the class of standard stationary
spacetimes. We begin by describing the central result in
\cite{LiNiren}. Let $(S,F)$ be a  complete Finsler $n$-manifold,
and $\Omega \subset S$ an open connected subset such that its
boundary $\partial \Omega$ satisfies  the H\"older condition
$C_{\rm loc}^{2,1}$. Let $G$ be the subset of $\Omega$ containing
the points where the closest point from $\partial\Omega$ is unique.
Then $\Sigma=\Omega\setminus G$ is  a subset of  the set of
points where the inner ``normal'' geodesics from $\partial \Omega$
do not minimize anymore (i.e., the cut locus). Now denote  by
$\ell(y)$ the length of such a inner normal geodesic from
$y\in\partial\Omega$ to the first hit in the cut locus $m(y)\in\Sigma$. In
\cite[Th. 1.5-Cor. 1.6]{LiNiren} the authors proved the following
(optimal) result:
\begin{thm}\label{LiNi} Under the ambient hypotheses above, for any $N>0$ the function
\[\partial \Omega\ni y\mapsto\min(N,\ell(y))\in\R^+\]
is Lipschitz-continuous on any compact subset of $\partial
\Omega$. As a consequence
\[{\mathfrak h}^{n-1}(\Sigma\cap
B)<+\infty\] for any bounded  subset $B$, where ${\mathfrak
h}^{n-1}$ denotes the $(n-1)$-dimensional Hausdorff measure.
\end{thm}
Now, let $A$ be a closed achronal hypersurface with boundary of  a spacetime
$(M,g)$ of dimension $(n+1)$. It is known that  any point $p$ in
$H^+(A)$ admits a {\it generator}, i.e., a lightlike geodesic
through $p$ entirely contained in $H^+(A)$ which is either
past-inextensible or has a past endpoint in the boundary of $A$
(see for example \cite[Prop. 6.5.3]{HawEll73}). Let us denote by
$N(p)$ the number of generators  through $p\in H^+(A)\setminus A$,
and $H_{\rm{mul}}^+(A)$ the {\em crease set} \cite{BK98, ChrGal}
i.e., the set of points $p\in H^+(A)\setminus A$ with $N(p)>1$. It
is known that $H^+(A)\setminus A$ is a topological hypersurface
which satisfies a Lipschitz condition \cite[Prop.
6.3.1]{HawEll73}; therefore, its non-differentiable points
constitute a set of zero ${\mathfrak h}^{n}$-measure (even though
this set may be highly non-negligible \cite{ChrGal}). Moreover,
the set of points where $H^+(A)\setminus A$ is not differentiable
 coincides with the crease $H_{\rm{mul}}^+(A)$ (see
\cite{BK98}). Using Th. \ref{LiNi} and Prop. \ref{p4.5}, the
following more accurate estimate
on the measure of this set is obtained. 
\begin{thm}\label{nondiffmeasure}
Let $(\R\times S,g)$ be  a standard stationary $(n+1)$-dimensional
spacetime with $S_0=\{0\}\times S$ Cauchy, and let $\Omega\subset S$ be an open
connected subset with $C^{2,1}_{{\rm loc}}$ boundary $\partial
\Omega$. Put $A= \bar \Omega$, $A_{t_0}=\{t_0\}\times A$, and let
$H_{\rm{mul}}^+(A_{t_0})$ the crease set of  $H^+(A_{t_0})$. Then, for
any compact (or Fermat bounded) subset $B \subset S$, we have that
\[{\mathfrak h}^{n-1}((B\times \R)\cap H_{\rm{mul}}^+(A_{t_0})
)<+\infty .
\]
\end{thm}

\section{Causality applied to Randers metrics}\label{randersfermat}
Consider now any Randers manifold $ (S,R)$  as
 defined in (\ref{randers}). In \cite{BCS04}, the authors use the
expression of a Randers metric as a Zermelo metric in order to
classify Randers metrics of constant flag curvature. Here we will
study Randers metrics with compact symmetrized closed  balls by exploiting
their expression   as Fermat metrics for a standard stationary
spacetime as in \eqref{standard}. Concretely:
\begin{equation}\label{randfermat}
\begin{cases}
 \beta=1,\\
 g_0=h-\omega\otimes\omega,\\ 
\end{cases}
\end{equation}
(see also \cite{BiJa08} for a
description of the equivalence between Randers, Zermelo and Fermat
metrics).

\subsection{Geodesic connectedness in Randers metrics.}
The first consequence, for Randers metrics, of the interplay with
Fermat ones is the following.

\begin{lemma}\label{lgeodesicconnection}
Let $(M,R)$ be a (connected) Randers manifold.
If its  symmetrized closed balls
are compact then it is convex.
\end{lemma}
\begin{proof}
Consider the expression of the Randers metric as a Fermat metric
described in \eqref{randfermat}. We know from Th.
\ref{causalstructure} $(\rm b)$ that the associated standard
stationary spacetime is globally hyperbolic and then causally
simple. By part $(\rm a)$ of the same theorem, the Fermat metric
is convex, i. e. there exists a minimal geodesic between every two
points $p,q\in M$.
\end{proof}
This lemma is interesting because its proof provides a geometrical
understanding of the compactness assumption on the symmetrized closed balls and it suggests the optimal
Finslerian result. Indeed,
Lemma \ref{lgeodesicconnection} can  be generalized  to
every non-reversible Finsler metric by using analytical techniques.
The key point is to prove that,  if the
symmetrized closed balls are compact,
the energy functional \eqref{E} of a Finsler metric satisfies the
Palais-Smale condition on the manifold $\Omega_{p,q}$ of
$H^1$-curves joining two given points $p$ and $q$ of $M$. By using
variational arguments, in \cite[Th. 3.1]{cmm}, it is proved that
if $(M,F)$ is forward or backward complete then the energy $E$ satisfies the
Palais-Smale condition on $\Omega_{p,q}$. But the forward or
backward completeness is only used  to show that, given a
Palais-Smale sequence, there exists a uniformly convergent
subsequence. To that end, it is enough to prove that the supports
of the sequence  are contained in a compact subset and, then, to
apply the Ascoli-Arzel\`a theorem. This follows from the
Finslerian Hopf-Rinow theorem,  using that the forward (or
backward)  closed balls are compact. But we can show easily that
the Palais-Smale sequence is contained in the intersection of two
closed balls $\bar{B}^+(p,r_1)\cap\bar{B}^-(q,r_2)$ for some
$r_1,r_2\in\R$. Therefore, by Prop. \ref{dscomplete}, it is enough
to assume that the symmetrized closed balls are compact. 
Once the Palais-Smale condition is satisfied,
a standard minimization argument based on the Deformation Lemma (see for
instance \cite{MawWil89}) applies, giving the existence of a geodesic connecting $p$ and $q$ and with length 
equal to $d(p,q)$. 
Moreover, by using Ljusternik-Schnirelmann theory, we 
obtain also the existence of infinitely many connecting geodesics
if $M$ is not contractible. Summing up:
\begin{thm}\label{tgeodesicconnection}
Any Finsler metric
with compact symmetrized
closed balls
is convex, i.e., for any  $p,q\in M$ there exists a geodesic
joining $p$ to $q$ with length equal to the distance $d(p,q)$.
Moreover, if $M$ is not contractible then   infinitely many
connecting geodesics with divergent lengths exist.
\end{thm}
Further developments of this result for manifolds with
boundary have been obtained in
\cite{BCGS}.
\begin{rem}
 The Finslerian Theorem of Bonnet-Myers (see for example
\cite[Theorem 7.7]{BCS01}) can be also extended by assuming
compactness of symmetrized closed balls rather than forward or
backward completeness. The proof of the theorem under this
hypothesis can be accomplished by following the same steps as in
\cite[Theorem 7.7]{BCS01}. As a consequence, this condition can
also be considered in Synge's Theorem (see \cite[Theorem
8.8.1]{BCS01}) and, then, in theorems which may be formulated by
using it implicitly, as
the sphere theorem (see Rademacher's version  focused on
non-reversible Finsler
metrics, 
\cite{Rad04b}).
\end{rem}
\subsection{Completeness of the symmetrized distance.}
The role of  compact symmetrized closed balls in previous results
 suggests to discuss when a Randers metric has complete
symmetrized distance. Along the way, we will obtain some necessary
conditions for a standard stationary spacetime to be globally
hyperbolic.
\begin{prop}\label{hcomplete}
Let $(M,R)$ be a Randers metric as in \eqref{randers} with complete symmetrized distance. Then the Riemannian manifold $(M,h)$ must be complete.
\end{prop}
\begin{proof}
It is enough to prove  that $\ds(p,q)\leq \dist_h(p,q)$, which can
be done as follows. Let $\ell_h, \ell_R$ denote, resp., the length
measured with $h$ and $R$. For any smooth curve $\alpha$ joining
$p$ and $q$ and its reverse curve $\tilde{\alpha}$, one has
$\ds(p,q)\leq \frac
12(\ell_R(\alpha)+\ell_R(\tilde{\alpha}))=\ell_h(\alpha)$. As for
every $\varepsilon>0$, $\alpha$ can be chosen such that
$\ell_h(\alpha)<\dist_h(p,q)+\varepsilon$, the required inequality
follows.
\end{proof}
Completeness of the Riemannian metric $h$ is  only a necessary
condition for the completeness of $\ds$, as the following
counterexample shows.
\begin{exe}\label{nosymcomplete}
Consider $\R^2$ with the Euclidean metric $\langle
\cdot,\cdot\rangle$ and the  sequence of points
$\{p_n=(0,n)\}_{n\in\N}$. For each $n\in\N$, choose a unit-speed
injective curve $\gamma_n=(x_n,y_n):[0,2]\rightarrow \R^2$ from
$p_n$ to $p_{n+1}$, with $\gamma_n|_{(0,2)}$ contained in the set
$0<x,\ n<y<n+1$. Consider also the $y$-symmetric curves
$\tilde{\gamma}_n=(-x_n,y_n)$. Let  $\varepsilon_n>0$
 small enough and $0<\alpha_n<1$, close to $1$, such that
$\varepsilon_n\!+1-\alpha_n\!<2^{-n-1}$. Choose functions
$\mu_n:[0,2]\rightarrow \R$, $\mu_n= \alpha_n \hat\mu_n$, where
$\hat\mu_n$ is a bump function equal to $1$ in a neighborhood of
$[\varepsilon_n,2-\varepsilon_n]$, and equal to $0$ in a
neighborhood of  $0$ and $2$. Finally, define a $1$-form  along the
curves $\gamma_n$, $\tilde{\gamma}_n$ as
$\omega_{\gamma_n(s)}(v)=-\langle\mu_n(s)\dot\gamma_n(s),v\rangle$
and $\omega_{\tilde{\gamma}_n(s)}(v)=\langle\mu_n(s)\dot
{\tilde{\gamma}}_n(s),v\rangle$, and extend it to a  $1$-form
$\omega$ on all $\R^2$ with norm strictly less than $1$ at every
point.

Now consider the Randers metric $R$ in  $\R^2$ associated to the
Euclidean metric and the $1$-form $\omega$. The
(non-converging) sequence $\{p_n\}_n$ is Cauchy for the
symmetrized distance $\ds$ of $R$. In fact, both $\dist(p_n,p_{n+1})$
and $\dist(p_{n+1},p_n)$ are smaller than $2^{-n}$, as can be checked
by computing the Randers length of the curves $\gamma_n$ and
$t\mapsto\tilde\gamma_n(2-t)$. For example, for the first one:
\begin{multline*}\ell_R(\gamma_n)=\int_0^2
\left(\omega\big(\dot\gamma_n(s)\big)+\sqrt{\langle\dot\gamma_n(s),\dot\gamma_n(s)\rangle}\right)\df s\\
\leq
2\varepsilon_n+\int_{\varepsilon_n}^{2-\varepsilon_n}
\left(-\alpha_n\langle\dot\gamma_n(s),
\dot\gamma_n(s)\rangle+\sqrt{\langle\dot\gamma_n(s),\dot\gamma_n(s)\rangle}\right)\df
s\\< 2\varepsilon_n+ 2(1-\alpha_n)<2^{-n},
\end{multline*}
as required.
\end{exe}
Finally, as a straightforward consequence of
Th. \ref{causalstructure} $(\rm b)$ and  Propositions
\ref{dscomplete} and \ref{hcomplete}, we have:
\begin{cor}\label{globhyp}
Let $ (\R\times S,\g)$ be a globally hyperbolic standard stationary
spacetime as in \eqref{standard}. Then  the Riemannian metric on $S$
$$ h=\frac{1}{\beta}g_0+\frac{1}{\beta^2}\omega\otimes\omega,$$
is complete.
\end{cor}

\begin{rem}  In the static case ($\omega=0$), the converse is
also true. Nevertheless, this is not true in general  (a
counterexample would follow easily from Example~\ref{nosymcomplete}).
\end{rem}
\subsection{Randers metrics with the same pregeodesics. } \label{s5.3}
Recall first that  different Randers metrics with the same
pregeodesics can be obtained by adding an exact $1$-form $\df f$ with
small enough norm, concretely, such that
$\df f(v)<1$ for all $ v\in TS$ with $ R(v)=1$.
We denote  by  $R^f$ the Randers  metric $R^f(v)=R(v)-\df
f$. In order to check that the
pregeodesics of $R$ and $R^f$ coincide, notice that those joining two fixed points
$p,q\in M$ are the critical points of the Finslerian length. As
for any curve $\alpha:[a,b]\rightarrow M$ joining $p$ and $q$  we
have that
\begin{equation*}
\int_a^b R^f(\dot\alpha(s))\df s =\int_a^b
R(\dot\alpha(s))\df s+f(p)-f(q),
\end{equation*}
the critical curves of $R$ and $R^f$ coincide.
 Moreover,  the symmetrized distance
associated to $R^f$ coincides trivially with the one associated to
$R$.\footnote{This fact can be used to find an example of a
Randers metric with complete symmetrized distance but with
non-compact symmetrized closed balls as in Example
\ref{counterhopfrinow}, in a more  elegant way.}

From the viewpoint of SRC in the Introduction, the Randers
metric
$R$ can be seen as a Fermat one $ F_g $ for a spacetime $(M=\R\times
S, \g)$, and the function $f: S\rightarrow \R$ yields a section
$S^f=\{(f(x),x): x\in S\}\subset \R\times S$.  As we show
in the Introduction, if this section is spacelike then it induces
a different Fermat metric  on $S$ (as well as  a different
splitting of $(M,g)$ as a standard stationary spacetime with fixed
$K$). The next  results show that $S^f$ is spacelike iff $R^f=R-\df
f$ is Randers  and the corresponding Fermat metric $ F_{g^f}$ on
$S$ coincides with $R^f$.

\begin{prop}\label{multfermat0}
 Let $(S, R)$ be a Randers manifold and $(\R\times S, g)$ be the standard stationary spacetime associated to it via \eqref{randfermat}.
Let $f:S\rightarrow \R$ be a smooth function, then
$S^f=\{(f(x),x)\in\R\times S\,:\,x\in S\}$ is a spacelike
hypersurface if and only if
\begin{equation}\label{supremo}
 \sup_{v\in TS,\ R(v)=1}\df f(v)<1
\end{equation}
and, in this case, $R^f=R-\df f$ is also a Randers metric on $S$.
\end{prop}
\begin{proof}
 Let $\tilde f\colon S\rightarrow S^f\subset \mathbb{R}\times S$
be the map defined as $\tilde f(x)=(f(x),x)$. The tangent space at $\tilde f(x)$ to
$S^f$ is given by $T_{\tilde{f}(x)}S^f=\{\tilde{\xi}=(\df
f_x(\xi),\xi) : \xi\in T_xS\}$.  Evaluating $g$ on vector fields $\tilde\xi\in TS^f$,  we get
\begin{equation}\label{gprimozero}
 g(\tilde{\xi},\tilde{\xi})=g_0(\xi,\xi)+2\omega(\xi)\df f(\xi)-\df f(\xi)^2.
\end{equation}
By definition, $S^f$ is spacelike   if and only if
$$g_0(\xi,\xi)+2\omega(\xi)\df f(\xi)-\df f(\xi)^2>0,$$
for every $\xi\in T S,\  \xi \neq 0$. Since $g_0=h-\omega\otimes\omega$, this  condition is
equivalent to
\[
-R(-\xi)<\df f(\xi)<R(\xi),
\]
and then to
$$\df f(\xi) < R(\xi),$$
for every $\xi\in TS,\ \xi \neq 0$.
\end{proof}
\begin{prop}\label{multfermat}
 Let $(S, R)$ be a Randers manifold and $(\R\times S, g)$ be the standard stationary spacetime associated to it via \eqref{randfermat}.
Let $f:S\rightarrow \R$ be a smooth function and assume that
$S^f=\{(f(x),x)\in\R\times S\,:\,x\in S\}$ is a spacelike
hypersurface. Let  $g'_0$ be the Riemannian metric induced by $g$ on $S^f$ and $\omega'$ be the $1$-form on $S^f$ defined as $\omega'(\tilde \xi)=g(K,\tilde\xi)$, for all $\tilde \xi\in TS^f$. Consider the standard stationary metric $g^f$ on $\R\times S$ defined by $g^f_0=\tilde f^* g'_0$, $\omega^f=\tilde f^* \omega'$, $\beta=1$,  where $\tilde f\colon S\rightarrow S^f$ is the map $\tilde f(x)=(f(x),x)$. Then the Fermat metric $F_{g^f}$ on $S$ associated to $g^f$ is equal to  $R^f=R-df$.
\end{prop}
\begin{proof}
 For any $\xi\in TS$ we have
\[g_0^f(\xi,\xi)=g'_0(\df \tilde f(\xi),\df\tilde f(\xi))=g_0(\xi,\xi)+2\omega(\xi)\df f(\xi)-\df f(\xi)^2.\]
Moreover
\[\omega^f(\xi)=\omega'(\df \tilde f(\xi))=g(K,(\df f(\xi),\xi))=
 \omega(\xi)-\df f(\xi)
\]
The Fermat metric associated to $g^f$ is equal to
\begin{align*}
F_{g^f}(\xi)&=\omega^f(\xi)+\sqrt{g_0^f(\xi,\xi)+\omega^f(\xi)^2}\\
&=\omega(\xi)-\df f(\xi)+\sqrt{g_0(\xi,\xi)+\omega(\xi)^2}\\
&=\omega(\xi)-\df f(\xi)+\sqrt{h(\xi,\xi)}=R^f(\xi).
\end{align*}
\end{proof}
\begin{thm}\label{randerscompleteness}
Let $(S,R)$ be a Randers metric. There exists $f:S\rightarrow \R$
such that the Randers metric given by
\begin{equation}\label{randersf}
R^f(v)=R(v)-\df f(v)
\end{equation}
for $v\in TS$ is geodesically complete if and only if the
symmetrized closed balls of $(S,R)$ are compact.
\end{thm}
\begin{proof}
Consider the standard stationary spacetime $(\R\times S,  \g)$
associated to the Randers metric $R$ as described in \eqref{standard}
and \eqref{randfermat}. By  Th. \ref{causalstructure} (b), 
compactness of the symmetrized closed balls is equivalent to
global hyperbolicity. This property is also equivalent to the
existence of a smooth spacelike Cauchy hypersurface
\cite{BeSan05}, which can be  given as the graph of
some smooth function $f:S\rightarrow\R$. Moreover,  recall that a spacelike
hypersurface obtained as a graph $S^f$ is Cauchy iff the Fermat
metric $ F^f$ associated  to it is forward and backward complete (Th.
\ref{th2}). This is equivalent to the completeness of the pullback  of $F^f$
on $S$ through the map $\tilde f$.  From  Prop. \ref{multfermat},  the pullback metric is equal to $R^f$.
\end{proof}
\subsection{Cut loci of Randers metrics via Cauchy horizons}
In Subsection \ref{s4.3},
 the properties of the Fermat distance
 from some subset $A$ yield consequences on the
 horizon corresponding to $A$. Next, we will see that the
 correspondence is also fruitful in the converse way. In fact,
the applications of general results on Cauchy horizons for
Riemannian Geometry were already pointed out in \cite{CFGH}. Here,
this will be extended to Finsler Geometry.

 Let $(S,R)$ be a connected Randers manifold, not necessarily forward or backward
complete. Given any closed subset $C\subset S$,  the distance
function $\rho_C\colon S\rightarrow [0,+\infty)$ is the  infimum
of the lengths of the smooth curves in $S$ from\footnote{ All
 the results will be
 obtained for the distance $\rho_C$ from $C\subset M$ to a point $p\in
 M$. Analogous results
 hold for the distance from a point $p\in M$ to the subset $C\subset
 M$ --they are reduced to the former case by considering the reverse Finsler metric.}
$C$ to $p$. The function $\rho_C$ is Lipschitz, more precisely
$|\rho_C(p)-\rho_C(q)|\leq 2\,\ds(p,q)$. Let $I\subset [0,+\infty)$
 be a (non-empty) interval. We say that $\gamma:I\rightarrow S$ is a
{\em $C$-minimizing segment} if it is a unit  speed geodesic such
that $\rho_C(\gamma(s))=s$ for all $s\in I$. We emphasize that the
interval $I$ (which may be open, half open or closed) may not
contain $0$. Reasoning in the Finsler case as  in
\cite[Prop. 9]{CFGH} for the Riemannian one, we have:
\begin{prop}\label{existenceC}
 Every $p\in S\setminus C$ lies  on at least one $C$-minimizing segment.
\end{prop}
From now on we will assume that $C$-minimizing segments are
defined in their maximal domain. We say that a $C$-minimizing
segment has a {\it cut point} iff its interval of definition is of
the form $[a,b]$ or $(a,b]$ with $b<+\infty$ being then
$p=\gamma(b)$ the cut point. The set of  all the cut points is
called the {\it cut locus} of $C$ in $S$, denoted  ${\rm Cut}_C$.
For any $p\in S\setminus C$ let $N_C(p)$ be the number of
$C$-minimizing segments passing through $p$. By Prop.
\ref{existenceC}, $N_C(p)\geq 1$ for every $p\in S\setminus C$,
and it is easy to see that if $N_C(p)\geq 2$ then $p\in {\rm
Cut}_C$.

Now (taking into account formulas \eqref{randfermat}),  consider
the standard stationary metric  constructed for the {\em reverse}
Randers metric $\tilde{R}$   so that the
 past-pointing lightlike geodesics correspond to the geodesics of
$(S,R)$. We will focus in the ``lower half'' spacetime
$(M=\big(-\infty,0)\times S, \g\big)$. These choices are
convenient because we will use the general notion of horizon in
\cite{CFGH}, i.e,  a {\em future horizon} $\mathcal H$
is\footnote{These requirements would be fulfilled by  the
 horizons of Cauchy developments in Section \ref{s4.3}, if
one removes some parts of the spacetime; for example, for $A$
closed, $H^-(A)\setminus A$ would be a future horizon of
$M\setminus A$.} an achronal, closed, future null geodesically
ruled topological hypersurface. Here {\em future null geodesically
ruled} means that each point $p\in \mathcal H$ belongs to a future
inextensible lightlike geodesic $\Gamma \subset\mathcal H$, i.e. a
{\em null generator} $\Gamma$ of $\mathcal H$.
 Let us call $\mathcal H$
the graph of $-\rho_C$ in $M$, that is
\[\mathcal H=\{(-\rho_C(x),x):x\in S\setminus C\},\]
which is a future horizon in the sense above. Up to
reparametrization, the null generators are precisely, the curves
$s\mapsto (-\rho_C(\gamma(s)),\gamma(s))=(s,\gamma(s))$, where $\gamma$ is a
$C$-minimizing segment  of $(S,R)$.  Then, the number $N_C(x)$ of
C-minimizing segments through $x$ coincides with the number
${N(-\rho_C(x),x)}$ of null generators of $\mathcal H$
 through  the point $(-\rho_C(x),x)$. In
addition, the set ${\mathcal H}_{\rm end}$ (given by the past
endpoints of the null generators of $\mathcal H$) coincides with
the set $\{(-\rho_C(p),p):p\in {\rm Cut}_C\}$.  After a result by
Beem and Kr\'olak (see \cite[Th. 3.5]{BK98} and also \cite[Prop.
3.4]{ChrGal}) a point $p\in\mathcal H$ is differentiable iff
$N(p)=1$. Then, as a consequence we obtain:
\begin{thm}\label{nondiff}
Let $(S,R)$ be a Randers manifold, and $C\subset S$ a closed
subset.  A point $p\in S\setminus C$ is a differentiable point of
the distance function $\rho_C$  from $C$ if and only if it is
crossed by exactly one minimizing segment, i.e., $N_C(p)=1$.
\end{thm}
As a final remark, notice that, in  the preceding discussion, the fact
that $\mathcal H$ is a Lipschitz hypersurface (and, then by Rademacher theorem, almost
everywhere differentiable) and the result in \cite[Th.
1]{CFGH} about the zero $\mathfrak{h}^n$-measure of the set of
smooth ends, yield  directly:
\begin{cor}\label{nullmeasure}
If $C$ is a closed set in an $n$-dimensional Randers manifold
$(S,R)$, then  $\mathfrak{h}^n({\rm Cut}_C)=0$.
\end{cor}
We observe that the result in \cite{LiNiren} (see Th.
\ref{LiNi} above) says that, when the subset $C$ is regular
enough, then the Hausdorff dimension of ${\rm Cut}_C$ is at most
$n-1$. As far as we are aware it is not known if there exists a
subset $C$ such that the Hausdorff dimension of ${\rm Cut}_C$ is
equal to $n$.
\section{Appendix: the symmetrized distance and its path metric space}
As we have commented in Section \ref{finslermetrics}, the
Finslerian distance  $\df_F$ of a non-reversible Finsler manifold
$(M,F)$ is not a true distance, as it is non-symmetric.   One can symmetrize 
it to obtain the so-called symmetrized distance in
(\ref{symdist}), but then the analog to Hopf-Rinow theorem does
not hold (see the counterexample \ref{counterhopfrinow}).  Indeed, $\ds$
is {\em not} constructed as length metric (see
\cite{gromov}). In fact, we can construct from the symmetrized
distance its associated length metric as follows. Given a
continuous curve $\alpha:I\rightarrow M$ with $I\subset \R$ an
interval, an arbitrary subset of $\R$, we define the {\it
dilatation} of $\alpha$, ${\rm dil}(\alpha)$ as
\[{\rm dil}(\alpha)=\sup_{s,t\in I;s\not=t}\frac{\ds(\alpha(s),\alpha(t))}{|s-t|},\]
and the {\it local dilatation at $t_0\in I$} as
\[{\rm dil}_{t_0}(\alpha)=\lim_{\varepsilon\rightarrow 0}{\rm dil}(\alpha|_{(t_0-\varepsilon, t_0+\varepsilon)}).\]
As $M$ is a  (connected) manifold, we can consider just the class
of piecewise smooth curves,  and  define the length associated to
$\ds$ on $[a,b] \subset I$ as:
\[ l(\alpha)=\int_a^b{\rm dil}_t(\alpha)\df t .\]
 This length determines a length metric $\ds^l$ and, then, the path-metric space $(M,\ds^l )$ for which
the Hopf-Rinow theorem does hold  (see \cite[pp. 2--9]{gromov}).
For a Randers metric, using that $\ds(p,q)\leq \dist_h(p,q)$ for
every $p,q\in M$ (see the proof of  Prop. \ref{hcomplete}), and the definition
of length distance, we can easily deduce that $\ds(p,q)\leq
\ds^l(p,q)\leq\dist_h(p,q)$ for every $p,q\in M$. We will show
that actually $\ds^l=\dist_h$.
\begin{prop}
Let $ R(v)=\sqrt{h(v,v)}+\omega(v)$ be a Randers metric and $\ds$ the
symmetrized distance of $d_R$. Then the length distance $\ds^l$
associated to $\ds$ coincides with the distance $\dist_h$.
\end{prop}
\begin{proof}
It is enough to prove that given a curve $\alpha$,
parametrized by the $h$-length,
the local
dilatations ${\rm dil}^s_{t_0}(\alpha)$, ${\rm
dil}^h_{t_0}(\alpha)$ for  $\ds$ and $\dist_h$ satisfy:
\begin{equation}\label{00}{\rm
dil}^s_{t_0}(\alpha) \geq {\rm dil}^h_{t_0}(\alpha) (=1),
\end{equation} as the converse follows from $\ds^l\leq \dist_h$.
Consider a convex neighborhood $U$ of $\alpha(t_0)$ in the
Randers metric $R(v)=\sqrt{h(v,v)}+\omega(v)$, $v\in TM$. We can assume that the closure of $U$ is compact and contained in a chart $(\tilde U, x)$
such that $x(\tilde U)$ is a Euclidean ball in $\R^n$. Moreover,
we can take a constant $C>0$ such that
$$\frac{1}{C}|\df x (v)|\leq \sqrt{h(v,v)}\leq C R(v) \quad \forall v\in T
\tilde U,$$
where 
$|\cdot |$ is the Euclidean norm in $\R^n$. Consider an interval $I$ such that $\alpha(I)\subset
U$. Given $s,t\in I$, let $\gamma_1$ be the Randers
pregeodesic in $U$ from $\alpha(s)$ to $\alpha(t)$  and
$\gamma_2$ the Randers pregeodesic in $U$  from
$\alpha(t)$ to $\alpha(s)$ both defined in $[0,1]$ and
parametrized with constant $h$-Riemannian speed. Let $\gamma$ be
the closed curve defined in $[0,1]$ as $\gamma_1$ and in $[1,2]$
as $\gamma(t)=\gamma_2(t-1)$. Then
\begin{align*}
\frac{\ds(\alpha(s),\alpha(t))}{|s-t|}&=\frac{1}{2|s-t|}\int_0^2
R(\dot\gamma(\mu))\df \mu\\ &\geq
\frac{\dist_h(\alpha(s),\alpha(t))}{|s-t|}+\frac{1}{2|s-t|}\int_0^2
\omega(\dot\gamma(\mu))\df \mu ,
\end{align*}
and (\ref{00}) will follow if
\begin{equation}\label{limitest}\lim_{s,t\rightarrow t_0} \frac{1}{2|s-t|}\int_0^2 \omega(\dot\gamma(\mu))\df \mu=0.\end{equation}
Write in the chosen coordinates $\omega= \sum_i \omega_i\df x^i$
and 
$x\circ \gamma
= (\gamma^1,\ldots,\gamma^n)$ so that
 $|\df x(\dot \gamma)|^2= \sum_i (\dot \gamma^i)^2$.  Recall that
$h(\dot\gamma(\mu),\dot\gamma(\mu))$ is constant in $[0,1]$ as
well as in $[1,2]$ and, so, for all  $\mu\in [0,2]$:
\begin{multline}\label{precedent}
|\df x(\dot\gamma(\mu))|\leq C
\sqrt{h(\dot\gamma(\mu),\dot\gamma(\mu))}
< C\,\ell_h(\gamma) \leq C^2 \, \ell_R(\gamma) \\
=2C^2\, \ds(\alpha(s),\alpha(t))\leq 2
C^2\,\dist_h(\alpha(s),\alpha(t))\leq 2 C^2|s-t|,
\end{multline}
where $\ell_h$ and $\ell_R$ are the lengths associated to $h$ and
$R$ respectively. 
So, the mean value theorem and \eqref{precedent} imply the
existence of $\mu_i\in (0,\mu)$ such that:
\begin{multline*}
\omega(\dot \gamma(\mu))= \sum_{i=1}^n \omega_i(\gamma(\mu))\dot
\gamma^i(\mu) \\ = \sum_{i=1}^n \omega_i(\gamma(0))\dot
\gamma^i(\mu) + \sum_{i,j=1}^n \frac{\partial\omega_i}{\partial
x^j}(\gamma(\mu_i)) \dot\gamma^j(\mu_i) \dot \gamma^i(\mu ) \mu\\
\leq \sum_{i=1}^n \omega_i(\gamma(0))\dot \gamma^i(\mu) + 2\tilde
C (s-t)^2 ,
\end{multline*}
for some constant $\tilde C$ independent of $t, s$. Integrating
this expression:
$$\frac{1}{2|s-t|}\left|\int_0^2\omega(\dot\gamma (\mu))\df
\mu\right| \leq \tilde C |s-t|,$$ and \eqref{limitest} follows, as
required.

\end{proof} Observe that Example~\ref{counterhopfrinow}
yields a Randers manifold $(S,R)$ where $\ds$ is complete with $S$
$\ds$-bounded, non-compact and $\dist_h$-unbounded. The involved
Hopf-Rinow type relations are summarized in the following
diagram:
\[\xymatrix{
(\R\times S,g)\,\,\text{Causally simple}\ar@{<->}[rrr]^{\text{Th. 4.3}}&&&(S,R)\,\,\text{convex}\\
(\R\times S,g)\,\,\text{Glob Hyp}\ar@{<->}[rrr]^{\text{Th. 4.3}}\ar[u]^{\text{Trivial}}&&& (HB)_s\,\text{holds}\ar[u]_{\text{Th. 5.2}}\ar[d]^{\text{Prop. 2.2}}\\
&&&\ds\,\text{complete}\ar@{->}[d]^{\text{Trivial}}\ar@{->}[dlll]_{\text{Trivial}}\\
\dist_h\,\text{complete}\ar@{<->}[rrr]^{\text{Prop. 6.1 ($\dist_h=\ds^l$)}}&&&\ds^l\,\text{complete}\\
&&&
}\]

\end{document}